# A unified interpretation of probability

Louis Vervoort, PhD

*29.09.2026*

*Higher School of Economics, Moscow*

***Summary.*** *In this article I revisit the frequency interpretation of probability of von Mises, and argue that it contains the key elements of any coherent interpretation of scientific probability – notwithstanding the heavy criticisms von Mises collected from mathematicians and philosophers alike. I propose a modified definition that is not based on von Mises' 'infinite collectives' but retains an essential ingredient of his interpretation, namely that probability can only quantitatively be defined for events that are, or can be, repeated in similar conditions and that exhibit frequency stabilisation. New is that the mentioned 'conditions' should be 'partitioned'. Thus, I will partition probabilistic systems into object and environment, or in object, initiating, and probing subsystem, and show that such partitioning solves a series of classic problems. Other consistent interpretations can be unified with frequentism, e.g. Bayesianism à la Jaynes. I explore in detail the tantalising subjective or operationalist touch of probability, and argue that it can be put on objective grounds. One further upshot of the model is that several quantum puzzles can be seen as deriving from an imprecise interpretation of probability.*

## 1. General introduction

Philosophers are since a long time interested in the meaning, or interpretation, of the concept of probability, as witnessed by the enormous literature that exists on the topic. When people with a scientific training hear about this 'problem of interpretation', their first reaction is often one of surprise (as was my case). *Why is there a problem?* Are all problems of probability not adequately solved by the well-understood calculus of probability, based on Kolmogorov's axioms? That there *is* an interesting problem-space can be recognised by following two considerations – to which usually not much attention is paid in science curricula (regrettably, I now believe). First, in textbooks of applied probability, the science student learns that two or more interpretations of probability exist, notably the 'classic' interpretation of Laplace, applicable to chance games (using dice, cards, urns with balls, etc.), and the ‘frequency’ interpretation, applicable to basically anything else. Especially modern textbooks add other interpretations, notably the subjective interpretation. Rough definitions of these are the following. According to the classic Laplacian view, the probability of event E occurring in a probabilistic experiment is the ratio $N_E/N$, where $N_E$ is the number of outcomes in which E occurs, and N the total number of different outcomes. (This is a rough definition; it presupposes for instance that the outcomes are equally likely and mutually exclusive.) According to the frequency interpretation, the probability of E is the relative frequency with which E arises in long (ideally infinite) series of repeated experiments. Both these interpretations are ‘objective’. According to the ‘subjective’ or Bayesian interpretation of probability, probability values are subjective measures of belief; different persons may attribute different probabilities to the same event. Several other interpretations exist, and each in several variants. As a consequence, a philosophically-inclined student, consciously or unconsciously driven by an urge for unification, might wonder: *Why isn’t there one overarching interpretation?*

I will argue, following von Mises, that probability can be seen not only as a mathematical theory, but also as a physical one. In a sense, probability theory is the most encompassing physical theory we have. All physical systems comply with it: obviously all probabilistic systems comply with it[1]; and deterministic systems are a special case of probabilistic systems characterised by probabilities 0 or 1. This fact already seems to suggest that there should be *one* fundamental interpretation of (scientific) probability *as a basic physical property or category*. But can one really assert this intuition in view of the plethora of counterarguments and pluralistic positions that thrive in the literature? Isn’t that a thought too bold to entertain? Possibly, but according to the MCT-thesis we should take it seriously.

A second eye-opening consideration as to whether there is a real problem, stems from the intriguing chapter of the paradoxes of probability theory. Applied probability is rife with such paradoxes: problems the solutions of which were disputed even among famous mathematicians, including the fathers of probability theory. Many mathematics books on probability contain at least a few: Bertrand’s paradox, the Borel–Kolmogorov paradox, the Monty Hall problem, etc. Together with the first consideration, these paradoxes show that the *application* of measure-theoretic probability theory (Kolmogorov's axiomatic theory) to *real-world* problems is non-trivial. But it *can* of course be applied to real-world problems, for instance physical problems, say die throwing, chance games, or countless problems of statistical physics. In sum, these two considerations show that the interpretation of probability – What is probability beyond the mathematics? How to apply the theory? – is a genuine problem at the interface of mathematics, physics and philosophy. I now believe it is a particularly fascinating one, and a uniquely subtle one to boot. People who have delved into this question will not need to be convinced. Novices may now pause, and try to come up with their own answer to the above questions.

[1] I mean that the probabilities predicted by probabilistic physical theories, say quantum mechanics, obviously comply with probability theory.

Following the heuristic rule of Chapter 1, I believe the approach to tackle the problem of interpreting probability should be one of unification, as always: one should try to unify existing interpretations (as far as it is possible), and show that the unified theory can solve more problems, including the mentioned paradoxes and the very many problems discussed in the literature. I found it instrumental to guide this quest by following additional questions: What is a probabilistic *system* (in the physical sense)?, and: How does one extract probability values from empirical data? These questions put the search on firm empirical grounds, which seems not a luxury. Indeed, I believe it is not enough emphasised in the philosophical literature that in science *any* experimental probability is determined, measured, verified *as a relative frequency*. Recalling the verificationist principle of meaning of the logical empiricists, *and* the credo of synthetic philosophy that science is the ultimate judge of philosophy (at least for problems that have such an immediate link with science), this seems to give a head start to the frequency interpretation. But again, the real question is: does the latter allow one to unify and concomitantly to answer the largest number of well-posed questions? That is what will be studied in this chapter.

Some people immediately object that theories such as quantum mechanics may predict probabilities as real numbers that are not obviously ratios or frequencies (think of the square of a modulus of a quantum wave function, a transition amplitude, etc.). But to *verify* these numbers one always determines relative frequencies. And as most empirically determined quantities (say the magnitude of speed, energy, magnetic field, etc.), these relative frequencies necessarily have an experimental error or precision; they lie within a numerical range $\Delta Q$. As emphasised in the next section, it is the key assumption of my von Mises-style model, *and one of the best corroborated facts of science*, that $\Delta Q$ can be made as small as practically needed by increasing the number of trials. In other words, the measured probability $Q = P \pm \Delta Q$, where P is the 'real' or theoretical probability value and ΔQ the experimental precision or uncertainty. Any real number P between 0 and 1 (as predicted say by quantum mechanics) can be written under the form P = Q ± ΔQ with Q rational; and ΔQ can be made arbitrarily small by increasing the trials. In sum, one can determine ratios that arbitrarily closely approximate[2] the theoretical P.

Several authors have scrutinised the frequency interpretation (for reviews on this and other interpretations, see e.g. the excellent Fine (1973), von Plato (1994), Gillies (2000), Khrennikov (2008)). K. Popper and B. van Fraassen were at some point proponents of the frequency interpretation – before becoming adepts of other positions. I will have a closer look at Popper's arguments in favour of the so-called propensity interpretation in Section 5. In his (1980, pp. 190 – 194), van Fraassen analyses how to link physical experiments to probability functions, and proposes following definition (1980, p. 194):

> "The probability of event A equals the relative frequency with which it would occur, were a suitably designed experiment performed often enough under suitable conditions."

My definition will be seen to be broadly in agreement with the above one. But it appears that the notion of probability has so many hidden dimensions that a more elaborate model is needed in order to address the numerous extant problems. To mention only two of these implicit notions: I will give in particular a further analysis of what the 'suitable conditions' in van Fraassen's definition are, and what 'often enough' might mean. This being said, it will be seen that van Fraassen's condensed phrasing can function as an excellent mnemonic device.

Besides philosophers, the fathers of probability theory, from Pascal, Fermat, De Moivre, Laplace to von Mises and Kolmogorov, were often interested in the emerging philosophical questions. The analytic definition of probability I will propose in Section 4, the core of this chapter, can be seen as a variant of von Mises' definition (1928/1981, 1964). When trying to construct my interpretational model, I found Kolmogorov's and his pupil Gnedenko's standard references on the mathematical theory highly instrumental[3] (Kolmogorov 1933/1956, Gnedenko 1967).

[2] Note that if, in this verification process, the predicted probabilities do consistently not correspond to the measured frequencies (within the experimental error), the theory is at risk of being rejected.

[3] I benefited in particular from Boris Gnedenko's book on probability calculus. Gnedenko did not only contribute to the calculus (e.g. by elaborations of the Central Limit Theorem and in statistics); his work reflects his keen interest in the foundational issues. Also Kolmogorov insisted on the need of an interpretation of probability (cf. his (1933/1956) p. 9). Interestingly enough, the axiomatisation of probability theory and the discovery of quantum mechanics overlapped in time. Von Mises but also Kolmogorov and Gnedenko followed the development of quantum mechanics with interest, surely realising that probability theory governs an ever wider part of the natural phenomena.

Actually, Kolmogorov adopts an interpretation of probability that is frequentist (1933/1956, p. 3, see below); he explicitly refers to von Mises as his source of inspiration.

A last preliminary remark concerns the use of the word 'probability' in daily contexts. Even if there may be other uses of the notion of probability in general philosophy, psychology, and certainly everyday discourse[4], and even if some of these may be coherent, I will focus here on the notion as used in natural science. I will find as usual initial inspiration in physics, but the aim is to encompass natural science as a whole.

The chapter is organised as follows. In Section 2, I will anticipate some of the modifications I will introduce to von Mises' interpretation. In Section 3 I will analyse some simple probabilistic systems in order to introduce the main concepts of my definition, which will be presented in final form in Section 4. Problem-solving and comparison with other interpretations are done in the Sections 4 and 5. Section 6 concludes and proposes lines of further research.

## 2. Introduction for the professional. Preliminary remarks on similarities and differences with von Mises' model

It was A. Kolmogorov who published the first and simplest axiomatic system for probability in 1933 (1933/1956). This theory is generally believed to cover all probabilistic systems of the natural, applied and social sciences. At the same time Kolmogorov's theory – the calculus of probabilistic or random events – does not define what a 'probabilistic / random event' is and does not provide an interpretation of probability. For a deeper understanding of what probability is, and for applying the calculus to real world situations, one needs to resort to the ideas of other fathers of probability theory, such as Laplace, Fermat, Venn, von Mises, and to the philosophers having investigated this matter (for reviews, see notably Fine (1973), von Plato (1994), Gillies (2000), Khrennikov (2008)).

In the broader physics community the most popular interpretation surely is the frequency model – especially the limiting frequency version, usually attributed to Richard von Mises (1928/1981, 1964). But in philosophy and the foundations of quantum mechanics other interpretations, in particular the subjective interpretation, are increasingly popular. Suffices to consult contemporary works on the interpretation of probability in quantum physics[5] (and, to a lesser extent, in statistical mechanics) to realise how vivid and wide-ranging the debate is, also in the physics community. Thus, the objective versus subjective controversy is far from being settled (see e.g. Gillies (2000) and for an overview in physics, Beisbart and Hartmann (2011), Ch. 1). I will try to come to grips with this remarkably subtle debate in the following sections.

[4] In daily conversations, for instance when discussing "the probability that person X will do Y" it is more appropriate to use the word 'likelihood' rather than 'probability'; the former concept is not necessarily related to the strict rules of probability calculus. Thus, in principle, it is preferable to reserve 'probability' for the concept that complies with the mathematical theory.

Now, the dramatic underrepresentation of the frequency interpretation à la von Mises in the contemporary philosophical debate is very striking. This state of affairs is partly due to the fact that naïve frequency interpretations have extensively been analysed and found wanting (cf. e.g. Fine (1973), von Plato (1994), Gillies (2000), Hájek (2009), (2023) and references therein). More importantly, von Mises' calculus, based on the notion of 'collective', acquired both among mathematicians and philosophers the reputation of being complex and possibly incoherent. During the last century several alleged mathematical flaws were discussed in publications, but today most commentators seem to agree that von Mises' calculus can well be made consistent if amended and interpreted adequately (e.g. Gillies (2000) p. 104; von Plato (1994) p. 197; Khrennikov (2008) p. 27). Still, the damage was done.

Indeed, I believe that these numerous criticisms have unduly discredited von Mises' work. I will go further: I believe that neglecting the key ideas that are explicit or implicit in this work simply amounts to being blind for the most fundamental properties of scientific probability (beyond the mathematics). One could say that where Kolmogorov describes probability, von Mises explains it[6]. Von Mises' theory not only recovers Kolmogorov's axioms, it also explains what probability values or measures 'are', where probability comes from – *and it is the only one to achieve this in detail*. The subtlety of it, and hence the many uncharitable criticisms, are, I submit, in the first place due to the fact that the theory does not have to deal with objects (or events) *per se*, as usual physical theories, but with *collections* of systems or events ("mass phenomena and repetitive events", in von Mises' words ((1928 / 1981) p. v)). I conjecture that one needs a genuine change of mind-set to come to grips with this particularity, as we will study in this chapter.

In *rough* essence von Mises interprets the probability $P(R_j)$ as $\lim_{n\to\infty} n(R_j)/n$, where $n(R_j)$ is the number of events or trials that have result $R_j$ in a series of n trials ($j \in \{1,\ldots,N\}$). This concept of *limiting relative frequency*, for the limit of n tending to infinity, is considered problematic by several authors (see e.g. Neapolitan (1992), who disagrees). This is a remarkable state of affairs. By definition of limit,

$$P(R_j) = \lim_{n\to\infty}\frac{n(R_j)}{n} \Leftrightarrow \forall\varepsilon>0, \exists\, \text{N: } \forall n>N, \left|\frac{n(R_j)}{n}-P(R_j)\right|<\varepsilon\,. \qquad (1)$$

[5] See e.g. Caves et al. (2002), (2007), and Bub (2007).
[6] I thank Richard Gill for this apt phrase.

Mathematicians as Richter (1978, p. 426) and others (Neapolitan 1992) complain that the last inequality in (1) may not be satisfied. But this is in overt contradiction with empirical evidence: an overwhelming amount of experimental data has shown that in long (of course not infinite) experimental series the values $n(R_j)/n$ do tend towards a fixed number *in a way that complies with (1) – at least under conditions that I will specify in detail in the next sections*. (Let me warn in advance: convincing oneself that (1) is correct may demand quite some attention. A clear graphical illustration of the limiting behaviour of typical random variables is given in Gillies (2000), p. 93, Fig. 5.2.) A first somewhat trivial confusion regarding (1) is related to the observation that surely for some chosen n and small ε the difference $\left|\frac{n(R_j)}{n}-P(R_j)\right|$ might be larger than ε. But experience has shown that one should be able to find a number N large enough such that for all n > N the difference becomes smaller than any ε. This well-corroborated scientific hypothesis of frequency stabilisation – the 'Ur-phenomenon' of probability theory as von Mises calls it – is well captured by (1) *if one understands it correctly*, as argued in detail in the next sections.

A less obvious and frequent confusion is linked to the following. Note first that the existing empirical evidence corroborates the initial hypothesis (namely of frequency stabilisation as a real-world phenomenon) as strongly as any hypothesis of the natural sciences can be confirmed: hypotheses of natural science are never confirmed with absolute certainty. Now, of course, in some series of trials aiming at determining probabilities, it may be that there are changes in the limiting value. For instance, it may be that after very many repeated measurements of a magnetic field strength B in a given point X, leading to a constant result of say 0.5 Tesla (which would allow one to infer that $P(B_X = 0.5) \approx 1$), the next measurements in the series unexpectedly jump to 0.7 Tesla. *This in itself would not be enough ground to reject the hypothesis that, in the conditions in which the first series of measurements were done, there is a fixed field $B_X \approx 0.5$ Tesla in X.* If the next trials constantly indicate that now $B_X \approx 0.7$ Tesla, physicists would look for causes of the jump. The conditions (environment) of the experiment must have changed, or the magnetometer must have started to malfunction; 0.7 Tesla is the field strength in the new (stable) experimental conditions – which arose due to some known or unknown systematic cause. If 0.7 Tesla is only measured for extremely rare trials in the series (with probability measure 0), this jump would be termed 'erratic'; it would be due to an unknown, unsystematic cause and it would *not* alter the conclusion $P(B_X = 0.5) \approx 1$. Similarly, if in a given experiment that is known to be genuinely probabilistic a limiting frequency could not be identified, one would not reject the law of frequency stabilisation, but look for causes to explain the erratic behaviour. We will look at all concepts involved here in more detail in the next sections.

Thus, at the very basis of the frequency interpretation lies the scientific hypothesis that if certain well-defined experiments (trials) would be repeated an infinity of times, certain ratio's (involving outcomes of these trials) would converge to a limit[7]. Some read this as a strongly counterfactual statement ("infinite series do not exist"); *but note that that is not the only possible reading*. Definition (1) can also be read in a very pragmatic and factual way: the longer one makes the trial series, the closer one comes to a certain number (in a way compatible with (1)) – as confirmed by countless experiments in all areas of natural science. I therefore emphatically disagree with authors as A. Hájek, who comments: "Note that at this point we have left empiricism behind. A modal element has been injected into frequentism with this invocation of a counterfactual; moreover, the counterfactual may involve a radical departure from the way things actually are, one that may even require the breaking of laws of nature" (Hájek 2023, Sect. 3.4). Much to the contrary: frequentism, if well defined, is likely the only interpretation that takes the laws of nature seriously. I will look in quite some detail into the allegedly counterfactual touch in von Mises-style theories in the next sections.

Indeed, this stabilisation hypothesis suggests one of the most fascinating questions of philosophy: *Why does frequency stabilisation happen? Why is frequency stabilisation a ubiquitous property of nature?* To my personal taste, these (and related) questions point to some of the most intriguing problems of science and philosophy. I will return to them in the Epilogue.

For the sake of completeness, I very succinctly analyse some of the criticisms of von Mises' theory (mainly related to his use of the notion of infinity) in Appendix 4.1. It should be remarked, though, that von Mises dealt in detail with many of them (see in particular the introductions of his 1928/1981 and 1964). If probability theory is a physical theory, then it is useful to remember that the concept of infinity is ubiquitous in theoretical physics, even if it is an idealisation in calculations – *a useful one when the formulas correctly predict what happens in reality*.

[7] More on the modal aspect of probability in (van Fraassen 1980).

Nevertheless, in the following I will argue that it is not necessary to define frequency stabilisation in terms of convergence in the limit $n \to \infty$. One can define 'probability-as-measure-with-precision-$\delta$' in a pragmatic way without any counterfactual notion of infinity (cf. Section 4). Although I have never encountered a final argument against von Mises' use of infinity, introducing the notion of precision has the advantage (besides of avoiding putative subtleties) of allowing to define probability in a more inclusive way, as it is used in cosmology, meteorology, and other areas.

My frequency interpretation will be seen to be much simpler than von Mises' theory (1928/1981, 1964). A *real-world* (or *physical*) theory (T) for probability is the union of a calculus (C) and an interpretation (I): $T = C \cup I$. The theory 'I' may contain e.g. the philosophical hypothesis that the theory represents things or events 'out there'; more generally, it contains rules stipulating to which objects to apply the mathematics, and how to do this. Von Mises proposed both an interpretation I and a calculus C – a physical calculus one might say; not an abstract, purely mathematical theory as Kolmogorov's. Both his 'I' and 'C' are based on the concept of collective, an in principle infinite series of experimental results. His calculus strikes however by its complexity, and nowadays probably no-one would consider using it on a regular basis (see e.g. his treatment of de Méré's problem, von Mises (1928/1981) p. 58ff). Therefore, I will not use von Mises' calculus of collectives; my calculus is Kolmogorov's; both lead to the same results.

In particular, I believe that it is not necessary to resort to von Mises' notion of 'invariance under place selection' to characterise randomness. Von Mises' 'randomness condition' is, in essence, the assumption or postulate that the limiting frequencies of a collective remain invariant in sub-sequences obtained under 'place selections', i.e. certain functions defined on the original collective (von Mises 1964). This characterisation triggered numerous mathematical questions and criticisms – to which von Mises responded. Now, even if von Mises' work on invariance under place selection is mathematically highly interesting (cf. Appendix 4.1), *I believe it is superfluous for its goal of defining randomness*. Indeed, it is enough to define randomness via the notion of frequency stabilisation. To the best of my knowledge, the most relevant argument is due to R. Gillies, who proved that the assumption of frequency stabilisation *without* the randomness condition suffices to imply Kolmogorov's axioms (Gillies (2000), p. 112). This clearly points to the conclusion that this condition is indeed redundant, and that frequency stabilisation is the only relevant criterion for probabilistic randomness (besides unpredictability). I will elaborate this idea in Sections 3 and 4.

Besides a simplification, my model is also intended as a clarification of the frequency interpretation. It helps in avoiding the paradoxes to which probability theory is so sensitive. I will show that when the application of probability theory becomes subtle, one gains in *partitioning* probabilistic systems into subsystems: namely test object; initiating; and probing subsystem (or environment, in some contexts of inquiry). This partitioning allows solving a series of classic problems (Section 5): it solves Bertrand's paradox; it offers an explanation of the 'subjective' temptation and of how it can be avoided; and it allows unifying the frequency interpretation with other interpretations. Finally, and perhaps surprisingly, my frequency model allows to re-interpret many puzzles of quantum philosophy as direct consequences of an inadequate interpretation of probability.

A last preliminary remark concerns the tantalising subjective or epistemic or operationalist dimension of probability theory, which we will encounter at different occasions throughout this chapter. I will only come to a final conclusion in the last section. The 'subjective versus objective' debate manifests itself already in the question of what the subject matter of probability theory is – to what exactly to apply it? According to von Mises probability theory treats "problems in which either the same event repeats itself again and again, or a great number of uniform elements are involved at the same time" (von Mises (1928 / 1981) p. 11). I will push his answer a little further and arrive at the conclusion that probability theory applies to a special type of random events or systems, which I will term 'p-random' or 'ρ-random' – hence 'ρ-systems'. Since I attribute probabilities to systems, my objectivist ambitions are clear. Indeed, my definition will be seen to be applicable to both 'chance games' and natural probabilistic phenomena: both happen in real-world systems[8]. (This is the first and simplest synthesis of this chapter.) In this objectivist perspective, I assume, likely as most scientists, that these natural random phenomena, such as quantum phenomena, spontaneously happen in nature *in agreement with the laws of probability*

[8] What I will term, following von Mises, chance games are games or procedures based on man-made probabilistic systems, such as dice, urns containing balls, card decks, roulettes, etc. 'Natural' probabilistic phenomena are ubiquitous in nature and can be found in diffusion, population dynamics, fluid dynamics, quantum mechanics and, actually, in any branch of natural science.

*theory – also when nobody looks*[9]. Of course, this assumption is in agreement with scientific realism, which I have to adopt in order to avoid the need to conjure up contrived answers to such questions as: "why can all people who use the same experimental protocol measure the same probability values for probabilistic events?" For the realist the answer seems straightforward and coherent with all of science: because probability is an objective, ontological category of nature; it characterizes *objective chance* (more in §4.2, **C2** below). (For a well-developed and different view, see Gillies (2000, Ch. 8) and, in the case of statistical mechanics, Uffink (2011).)

And yet I believe that this is not all that must be said about this unusually subtle matter. *There is an epistemic aspect to probability too*; more precisely: a somewhat complete physical/interpretational theory of probability will have to highlight epistemic aspects of the application of the theory. According to a typical definition "Objective interpretations of probability […] take probability to be a feature of the objective material world, which has nothing to do with human knowledge or belief. Clearly the frequency and propensity interpretations are objective" (Gillies (2008), p. 2). I agree, and yet I will argue in Section 4 that there is an epistemic dimension in the application of probability theory, in this precise sense: the ascription of whether a system or event is probabilistic (indeterministic) or deterministic, or neither of both, does depend on the knowledge state of the attributor. But let me reserve my final analysis to the conclusion (Section 6).

Let us now enter into the heart of the matter. Readers familiar with the frequency interpretation might go quickly over Section 3 (but halt to read Kolmogorov's interpretation and his reliance on the notion of 'condition') and focus on Section 4.

## 3. Examples. Frequency stabilisation and partitioned conditions

Let us have a look at a typical probabilistic system, namely a die. The key question to start with is: when or why is a die 'probabilistic', or 'random' (or rather the throwing events or outcomes)? Simply stating in non-anthropocentric terms what a die throw is, brings already to the fore a few essential notions. In physical terms, a 'die throwing event' or 'die throw' consists of the 3-dimensional movement of a (regular) die, that is 1) characterised by the time evolution of its centre of mass and its three Euler angles, 2) *caused* or *initiated* by an 'initiating / initialising system' (e.g. a randomising and throwing hand, or automat), and 3) *probed* by an 'observing / probing system' (e.g. a table) that allows to observe (in general 'measure') an outcome or result R (one up, two up, etc.). In the case of die throwing, it is easy to realise that if we want to use the die as it should, i.e. if we want to be able to observe or measure the usual probability for the different results of the throws, we have to repeat the experimental series (i.e. the throws) not only by using the same regular die, or similar regular dies, but also by applying the usual 'boundary conditions'. Here I obviously do not refer to the detailed initial conditions of each individual throw, which are unknown, but to the conditions of our random experiment *as they can be stated in an experimental protocol*, and thus be communicated to other experimenters ("throw roughly in this manner, probe on a table", etc.). Countless experiments have shown that if the protocol is well-defined and well-repeated, the same probabilities will be generated. In a slogan: same experiment, same probabilities. (If this notion of protocol rings operationalist bells to some readers, I can only agree. But, as said, I believe this aspect can be put on an objective basis.)

---

[9] It seems to me that one more convincing argument in favour of this claim resides in the following. One can measure averaged quantities (such as average kinetic energy) on ensembles of systems (say molecules in a gas), quantities that can also be calculated by integrating over probability distributions. The favourable comparison between measurement and theory – demonstrated many times by physicists – is indirect but convincing evidence in favour of the probabilistic assumptions made in the calculations, *and* of the claim that probabilities are objective.

These boundary conditions, then, are related to our throwing, the table, and the environment in general. *Irregular conditions in any one of these three elements may alter the probability distribution*. We can for instance not substantially alter our hand movement, say by putting the die systematically ace up on the table: this would change the probability distribution for the six outcomes from the usual (1/6, 1/6, 1/6, 1/6, 1/6, 1/6) to (1, 0, 0, 0, 0, 0). Nor can we put glue on the table, position it close to our hand, and gently overturn the die on the table while always starting ace up – one of the outcomes could occur more often than in a ratio of 1/6. Nor can we do the experiments in water instead of air: again one can imagine situations in which the environment alters the usual probabilities.

It will make sense to isolate the mentioned elements, and to consider a random event as involving a *random system* containing three subsystems, namely 1) the (random) *test object* itself (the die), 2) the *initiating system* (the throwing hand), and 3) the *probing* or *observing system* (the table, and let's include the human eye). I will call such a composed random system a 'ρ-system' or probabilistic system ('ρ' stands for the probability distribution generated by the system, see the general definition in **C3** below). Just as one can associate subsystems to the random event, one can associate (composed) *conditions* to it, so conditions under which the random experiment occurs: *initiating / initialising* and *probing conditions*.

Here is another example, from darts. According to the frequency interpretation, if we want to determine the probability that a given robot, equipped with a throwing arm, hits a given section of a dartboard (suppose there are 10 sections $S_1$ to $S_{10}$), we let it throw say 1000 times, and count relative frequencies: $n_1$ hits in $S_1$, so $P_1$ (the probability the robot hits $S_1$ in the given conditions) is approximately $n_1/1000$, etc. Again, if we want to define probabilities for the experiment, the robot, the darts *and* the board need to 'behave well' and operate in well-defined constant conditions. The probabilities are defined *relative to these conditions*, referring to the test object (the darts), the initiating system (the throwing robot), and the probing or observing system (the board, and we can again include some optical detection system). For instance, if during the experiment the robot's output force decreases due to a failing component in its circuit, $P_1$ might be lower than $n_1/1000$, the value obtained under normal conditions[10]. All this may look almost trivial; and yet it will prove useful to explicitly include these observations in the definition of probability.

At this point an essential remark is in order. Instead of 'initiating system' and 'probing system', it can be more appropriate to speak of 'environment', namely in the case of spontaneous or natural probabilistic events (versus artificial ones, as the outcomes of die throws and chance games created by human intervention). Many random events occur spontaneously, without any known cause. A spontaneously disintegrating nucleus has probabilistic properties, for instance the lifetime at which it disintegrates. Such properties are, on usual interpretations, not 'initiated' by a person or a robot, as in die throwing. Neither are they probed, except when subject to lab experiments. But the nuclei disintegrate spontaneously according to a probabilistic pattern that is numerically well-defined, *only if the environment is well-defined and stable*. Changing the environment, e.g. by irradiating the nuclei or by strongly heating them, may very well change the probability distribution in question. In other words, if we want to experimentally determine the half-life of the nuclei, i.e. the probability of disintegration, we have to put them in the lab in well-defined (initialising) conditions of temperature, pressure, etc. and measure their properties in well-defined (probing) conditions, which scrupulously correspond to their natural environment – if we want to know their probabilistic behaviour 'in nature'. So, also in the case of natural stochastic phenomena the initial and final conditions re-appear, be it under a different guise: as the conditions describing the environment.

[10] Notice there is some arbitrariness in where we draw the line between the different subsystems. But that doesn't play any role in the following.

By the above partitioning in subsystems we have clarified the concept of 'conditions'. This concept is not explicit in von Mises' notion of collective, but it is explicitly mentioned by Kolmogorov (1933/1956, p. 3-4), in the brief paragraph where he discusses how to apply his axiomatic theory to real-world *experiments* (sic). In this paragraph Kolmogorov explicitly acknowledges von Mises as his main source of inspiration (1933/1956, p. 3); and it is clear that Kolmogorov's interpretation comes very close to von Mises'. Here is the key passage:

> "There is assumed a complex of conditions, [C], which allows of any number of repetitions [....]. Under certain conditions, which we shall not discuss here, we may assume that to an event A which may or may not occur under conditions [C], is assigned a real number P(A) which has the following characteristics:
>
> (a) One can be practically certain that if the complex of conditions [C] is repeated a large number of times, n, then if m be the number of occurrences of event A, the ratio m/n will differ very slightly from P(A).
>
> (b) If P(A) is very small, one can be practically certain that when the conditions [C] are realized only once, the event A would not occur at all."

(Notice that 'conditions' is used in two different meanings. With the 'conditions that are not discussed here', Kolmogorov presumably refers to von Mises' conditions for applicability of probability theory, namely his conditions of frequency stabilisation and 'randomness'. I will argue below that the former suffices; it encompasses the latter.)

The notion of condition also appears in Gnedenko's (1967, p. 21):

> "On the basis of observation and experiment science arrives at the formulation of the natural laws that govern the phenomena it studies. The simplest and most widely used scheme of such laws is the following: *Whenever a certain set of conditions C is realized, the event A occurs*."

And a little further (1967, p. 21):

> "An event that may or may not occur when the set of conditions C is realized, is called *random*."

By partitioning, we characterised these conditions in some more detail than in the mentioned works.

Interestingly, this interpretational notion of condition can immediately be linked to the calculus – which shows how essential it is. Indeed, even unconditional probabilities (P(R)) can be regarded as conditional *on the conditions or circumstances of realisation*, as is implicit in Kolmogorov (1933/1956, p. 3) and explicitly stated by Gnedenko (1967, p. 67). In a scientific context these circumstances can be described as variables assuming a certain value (in that sense they correspond to events); in the simplest case such variables C have dichotomic values, i.e. $C = c \in \{yes, no\}$ or $\in \{\pm 1\}$, say. Thus, the natural generalisation in the model I will propose below is to replace, if helpful, P(R) by P(R|C) where C is the set or n-tuple of all relevant variables that describe the initiating, probing, and/or 'environmental' conditions in which the event R occurs. (Caution: C refers here to variables that are *constant* for all considered trials of the experiment; C are *not* random variables within the experiment considered. If C represents random variables then in general $P(R) \neq P(R|C)$, of course, and one cannot replace P(R) by P(R|C).) For instance, in the example of the disintegrating nucleus, the probability that the disintegration time $t = t_1$ can be formalized as $P(t = t_1 |T, p)$, where T and p fix the temperature and pressure[11] (in certain intervals). *Note that if these conditions are not made explicit as in P(R|C), then they have to be considered as implicit* (some authors use the notation '$P_C(R)$'). Again, at face value all this is quite straightforward.

With this in mind, it seems not difficult now to deduce the *conditions of randomness* of a physical event as a die throw – the criterion to which von Mises spent so much mathematical effort. Suppose we throw a die 20, 30, in general n times ($n >> 6$), and that we note how often the 6 results $R_j$ ($j = 1,\dots, 6$; $R_j = 1,\dots, 6$) occur in these n throws. We can then determine the relative frequencies of the results $R_j$, namely the 6 ratios $n_j/n$ where $n_j$ is the number of throws that have result $R_j$ and n the total number of throws. Now, according to the frequency interpretation (and, I submit, according to the practice of any scientist trained in experimental statistics), in order that the die throws can be termed 'probabilistic', it is a *necessary condition* that the relative frequencies of the $R_j$ (the ratios $n_j/n$) converge towards a constant number when n grows. For a regular die thrown in regular conditions it is an empirical matter of fact that the six ratios $n_j/n$ converge towards 1/6. *If such frequency stabilisation would not occur, one cannot speak of (scientific) probability*. If for instance the die would sublimate during the experiment in an asymmetric manner, thus gradually losing weight on one side, or if we modify the usual die throwing in a manner that does not lead to frequency stabilisation, we simply cannot attribute a probability to the outcomes. In such irregular conditions, if we repeat the throws, the above-mentioned ratios will not converge but erratically jump from one value to another without showing the typical limiting behaviour of probabilistic systems. In von Mises' words (1928/1981, p. 12):

> "It is essential for the theory of probability that experience has shown that in the game of dice, as in all the other mass phenomena which we have mentioned, the relative frequencies of certain attributes become more and more stable as the number of observations is increased."

[11] In principle it is more precise to write not $P(t = t_1 |T, p)$ but rather $P(t = t_1 |T = \tau, p = \pi)$ with $\tau$ and $\pi$ numerical values (or intervals); so to explicitly relate the probabilities to events instead of variables. This is understood in the widespread abbreviated notation I often use.

Thus, frequency stabilisation gets the status of a principle, or postulate, expressing something quite akin to a law of nature – and yet on the surface different from the usual laws, since now human agency (repeated observations) is explicitly referred to. Probability-theory-as-a-theory-of-the-universe, at least in von Mises' perspective, has to assume this principle, as backed-up by empirical evidence. Of course, any real-world theory has to assume postulates, laws,… (recall Chapter 2 on the problem of induction).

As already mentioned, frequency stabilisation and empirical probability[12] are characterised by a precision ('δ'), which will be integral part of the general definitions (**DEF1&2**) I propose below. The experimental series $n_j/n$ stabilise as a function of n towards the constant number P *within an experimental precision* $\delta$ (typically in an irregular, oscillating manner, but with on average ever smaller amplitudes, in agreement with (1)). This means that from some value n on, $n_j/n$ will lie in the interval [P-δ, P+δ]. Typically, this precision is of the order $1/\sqrt{n}$, in agreement with the Central Limit Theorem (Feller 1991; see examples in Gillies 2000, Ch. 5 and p. 94).

Since both these notions of precision and repeated observations are part of my definitions **DEF1&2** of probability, one could be worried that they have an obtrusive operational or anthropocentric or subjective dimension. Now, on closer inspection, one will see that *this is just an enhancement of a dimension that is present in any physical theory*. Take for instance Maxwell's electromagnetism. Maxwell's theory (as von Mises') is partly a mathematical theory: the electromagnetic fields satisfy equations, which are considered exact. These equations allow to extract, in principle, perfectly specified values for the fields, having $\delta = 0$. But the values measured in experiments only approximately correspond to the theoretical values, because of non-ideal *conditions*, such as non-zero background fields, or background charges that are not infinitely far away, etc. The precise value of the measured property depends on the specific conditions in which the experiment is done. So there is an operational / anthropocentric / subjective / epistemic dimension in the *application* of usual physical theories too; experimental conditions and precision are fundamental ingredients of the application of the theory. As a consequence, there is also a counterfactual dimension in the use of Maxwell's theory: the exact, 0-error value of the fields could only be measured if other background fields would exactly be equal to 0. But 'exactly 0' is a counterfactual notion in physics. Perhaps an even more striking example comes from the definition of 'electric field' in electromagnetism. Most textbooks on electromagnetism define the electric field (at an arbitrary point R) along these lines: the field is the electric force per unit charge that *would* be experienced by a tiny positive test charge if it *would* be placed at R. In sum, counterfactuality is ubiquitous and harmless in science, as in probability theory – a fact that seems not enough acknowledged (see e.g. Hájek 2023, Sect. 3.4). (More on this intriguing issue in **C4** in Section 4 and in Section 6.)

[12] Experimental or empirical probability is the object of interpretational / physical theories as the one we are interested in; to be compared to the purely mathematical probability of Kolmogorov's theory. Ultimately they should refer to the same thing. Analogously, physical properties usually are represented in physical theories as error-free variables with an exact value, which almost never coincides *exactly* with the experimental value.

It is just that in probability theory the operational, subjective dimension is a bit more prominent and confusing, apparently. More precisely, for probability the experimental conditions have a more decisive, and yet less easily graspable, influence on the experimentally determined values. I submit that this is ultimately due to the fact that this theory is all about partitioned systems and their repetitions, not objects *per se*, something we are not used to. What adds to the mystery is this: probability theory is, I believe, a more fundamental theory than usual physical theories. It is so general that it can be abstracted in an axiomatised mathematical theory. That is not the usual fate of physical theories (more in the Conclusion).

Interestingly enough, the idea that empirical probability 'needs' the notion of precision on a fundamental level, seems also suggested by Kolmogorov's insights condensed in point (a) in the above quote of his (1933/1956, p. 3-4).

It may be useful to further illustrate the ideas of partitioned conditions and frequency stabilisation by a realistic example from the physics lab, namely recent experimental findings in fluid mechanics and nonlinear physics, which have sparked considerable interest in the physics community in the last decades (e.g. Bush et al. 2024, Papatryfonos et al. 2024). A group of experimentalists have discovered that oil droplets can be made to hover over an oil film (Couder et al. 2005, Couder and Fort 2006, Eddi et al. 2011, Bush et al. 2024). To that end, the oil film is made to vertically vibrate using an external motor. If small oil droplets are deposited on such a vibrating film they sometimes begin to horizontally walk over the surface, for indefinite time; actually they bounce so fast on the film that they seem to hover[13]. However, this stable walking regime only occurs in well-defined experimental conditions, i.e. for precise values of the physical parameters of the system, essentially the frequency and amplitude of the external vibration, the size of the droplet, the geometry of the oil film and bath, and the viscosities of film and droplet. If these parameters are fine-tuned and lie within precise ranges of values, well-documented by the researchers, the droplets walk horizontally; outside these value ranges the movement becomes erratic and/or the droplet is captured by the film. In the walking regime it is possible to experimentally determine certain probabilities, notably the probabilities $P_R$ that the walking droplet or 'walker' passes through a certain space region R. To that end, the physicists have counted relative frequencies (cf. e.g. Couder and Fort. 2006, Fig. 2-3) much as I illustrated above in the case of a dartboard: they determine $n_R/n$, where $n_R$ is the number of trials in which the droplet passes through R, and n the total number of trials. Now, before determining numerical values of probabilities as $P_R$, the physicists have spent weeks, months, and possibly years to identify the exact conditions in which a stable probabilistic pattern (i.e. probabilities $P_R$) occurs. To *prove* that the system is probabilistic (and not erratic) there is only one way: namely to prove that frequency stabilisation occurs. When experimentalists state numerical values of probabilities in publications, they (should) have ascertained that such a frequency stabilisation has occurred in their experiments (if ratios / frequencies would be presented as 'probabilities' without anterior verification that these ratios lie in a stabilising regime, scientific deontology would be violated). In physics, this typically manifests itself by the fact that a probability histogram (e.g. $n_R/n$ as a function of discrete values of R) is better and better approximated by a smooth curve when n is increased (see e.g. Couder and Fort (2006), Fig. 2-3 where both histogram and curve are shown); often one publishes just the curve. This curve, in other words these limiting values of $n_R/n$, are supposed to be 'the' probabilities, or good approximations of it[14].

[13] These walking droplets exhibit moreover a behaviour that strikingly imitates quantum behaviour, including double-slit interference, tunneling and quantisation of angular momentum. Couder et al. have quite convincingly shown that the origin of such a quantum-like behaviour lies in the wave-front that accompanies the hopping droplets (the wave is created by the external vibration and by the back-reaction of the bouncing droplets on the oil film). Thus, these researchers could claim that these walking droplets are the first realisation of a "particle + wave" system (Couder et al. 2005, 2006). Here we do not focus on these intriguing quantum-like features, but on the classical probabilistic features, which were also extensively investigated by these researchers.

A key point highlighted by this hydrodynamic system is that frequency stabilisation, i.e. a probabilistic pattern, *only occurs under well-defined conditions, and by no means always*. Very small variations[15] in some variable (e.g. in the external vibration frequency, the droplet size, slight air movements, etc.) may lead to completely different probabilities, or to erratic behaviour without stable relative frequencies, or to absorbance of the droplet by the film. One has not necessarily control over these variations. Therefore, if a physicist who is novice in the dynamics of oil droplets on films – almost anyone – is asked to determine the probabilities $P_R$ by experiments, she may well come to the conclusion, even after weeks of experimentation, that the system is *not* probabilistic; or, if she is more cautious, that she did not identify conditions in which the system is probabilistic.

Before generalising these illustrations in definitions, in the above jargon probabilistic / random systems include deterministic systems. Indeed, in deterministic systems frequency stabilisation will also occur. If we move a die by placing it systematically – deterministically – ace up on the table, the six $P(R_j)$ will trivially converge to a value (1, 0, 0, 0, 0, 0). We could exclude deterministic cases from the definition of random events / systems by demanding that, besides exhibiting frequency stabilisation, the events / systems should be unpredictable. Thus 'genuinely' random events show frequency stabilisation *and* are unpredictable – or perhaps one could replace the latter by another predicate that qualifies the 'unstructured', disordered nature of randomness. (In a sense the unstructured-ness is only partial, since probabilistic systems are still subject to the mathematical structure of probability theory.) Deterministic systems could be called 'trivially probabilistic': their probabilities have the trivial values 0 or 1.

## 4. Generalisation, definitions

### *4.1. To what to apply probability theory?*

[14] For other examples of experimental determination of probability values, cf. (Papatryfonos et al. 2024) and (Nikolaev and Vervoort 2023).

[15] The extreme dependence of physical properties on certain variables is typical for nonlinear systems, such as those investigated by the teams of Couder and Bush. The oil droplet bouncing on the film is governed by a nonlinear equation of movement, e.g. due to the so-called 'viscous friction' between droplet and film.

If the essential feature qualifying physical events or systems as random-in-the-sense-of-probabilistic is frequency stabilisation, then following (neglected) question pops up: are all so-called random physical systems genuinely random in the sense of probabilistic? In the case of systems *studied in physics*, one usually declares that these are either probabilistic (indeterministic) or deterministic; so here random can be considered identical with non-trivially probabilistic. However, the above question becomes more interesting in the case of 'physical systems' in a larger sense, especially if the systems in question contain humans or are influenced by humans. In this case it is clear that *not* all so-called random systems are random-in-the-sense-of-probabilistic. *Frequency stabilisation is by no means guaranteed for an arbitrary disordered or spontaneous event or system*. Consider the number N of patients in a certain hospital division on a given day of the week, or the number N of cars parked in a street; and consider the probabilities that these variables assume a certain value ($P(N = N_1)$). The precise value of N varies unpredictably, randomly over time. Perhaps it is possible to define conditions (including e.g. the month of observation, etc.) for which the corresponding probabilities do stabilise when long series of data are gathered; but the notion of 'long series of data' becomes ambiguous here. Indeed, there may be *systematic causes* that change the limiting frequencies when these are identified in real experiments: a new vaccine may be developed that changes $P(N = N_1)$ over the duration of the experiment, or a new shop is opened in the street increasing $P(N = N_1)$ if $N_1$ is low. For (simple) physical events as a die throw one has a clear idea of what the conditions are that need to be repeated (one can repeat them 'as often as one wants'); but for complex systems, in particular human-induced ones, the conditions leading to frequency stabilisation may be so complex (or the set of conditions may be so large) that *they cannot be defined in an experimental protocol – they are not known in detail*. In such cases, experimental ratios $\frac{n_{N1}}{n}$ can be computed and are often still called probabilities; but this ratio is not a genuine probability: probability theory and statistics are not reliably applicable to it. It may make sense to call such numbers 'quasi-probabilities'[16].

As a realistic example, temperatures and rainy days in a given city look random, they even look statistically random, but they don't show frequency stabilisation over long periods, e.g. due to climate change. Meteorological probabilities may be seen as an interesting limiting case. If a frequency 'at infinity' cannot be determined, one may wonder whether the notion of probability is still sound here. I believe it is, *if one introduces the notion of precision in the definition,* as is done below. Indeed, contemporary meteorologists heavily use statistical talk: they assume that the conditions at a given time in a given place are sufficiently similar to those of which they have a large data set, allowing them to infer approximate probabilities. Since these data sets are very large, and climate conditions change rather slowly – at least in the previous decades! – *it is likely that the key requirement is satisfied, namely that frequency stabilisation occurs within a reasonable precision interval* $\delta$*.* What is reasonable is an (inter-)subjective matter, depending on the standards of precision of the field of application; it depends on which use one wishes to make of the inferred probability. In any case, one has here not the precision one can attain in more simple physical systems, for which the sampling size – of data corresponding to conditions that can be specified in an experimental protocol – is virtually infinite. But this problem can be handled if (empirical) probability is defined relative to a precision as is done in the definitions **DEF1** and **DEF2** below.

[16] Such quasi-probabilities may still have some practical value. For instance, hospital managers may use them to estimate the number of medical doctors to appoint in a given week. But this estimation is likely imprecise and unstable.

In the former section we saw a concrete example from fluid dynamics that illustrates the daunting practical difficulty for finding conditions of probabilistic stabilisation. Arguably a majority of physical systems are so complex and non-linear that they are extremely sensitive to variations in the parameters that describe them. As illustrated above, this leads often or 'normally' to the impossibility of determining conditions of frequency stabilisation. Interestingly, this is again a statement of pragmatic nature – one can argue that an omniscient being could know such conditions.

In conclusion, frequency stabilisation happens often, but certainly not always. (A particularly beautiful mathematical argument for the ubiquity of probabilistic randomness is provided by the Central Limit Theorem – sometimes called the unofficial queen of probability theory.) So not all randomness is the randomness that is the object of probability theory. This justifies the introduction of the concept of 'ρ-randomness', to be contrasted with 'randomness' and 'non-ρ-randomness': ρ-randomness, or probabilistic or structured randomness, is randomness (chance) that is characterised by frequency stabilisation 'for large n'. Before one applies probability calculus to an event, one conjectures or has evidence that it is ρ-random – a feature that can ultimately only be confirmed by experiment. I use the symbol 'ρ': in physics the familiar way to characterise a variable R as random, is to state that it has a 'probability density function', a parameter often symbolised by ρ (ρ = ρ(R)). In physics variables usually are continuous, so R can take an infinity of values (more on this point in **C5** below).

### *4.2. Definitions*

The stage is now finally set for introducing the two essential definitions of my philosophy of probability, **DEF1** and **DEF2**. In the comments **C1 – C6** I elaborate on the auxiliary terms and assumptions of these definitions, and mention the most immediate results to be inferred from them – so the most immediate justifications of the definitions. Von Mises starts by defining a collective as a series of random results and defines probability with respect to a given collective. But I have just argued that probability cannot be attributed to any random series, but only to ρ-random series or ρ-systems. Let us therefore first define the concept of frequency stabilisation: it is logically prior to that of probability.

In accordance with the discussions above, I define probability (with precision δ) relative to an ordered 5-tuple < C, T, $S_C$, {$R_j$}, δ> or, in a more condensed way, to the 4-tuple < C, T, {$R_j$}, δ> ($S_C$ can in principle be derived from T and C). Here C refers to the repeated experiment done on the test object T, more precisely to the partitioned experimental conditions mentioned above. $S_C$ is the larger system that consists of T in the conditions C (it is the ρ-system already mentioned). {$R_j$} represents the set of possible results/outcomes/events $R_j$ of the experiment (j = 1, …, N). (The $R_j$ can always be chosen as discrete numerical values or labels of continuous intervals of the real numbers.)

**DEF1 (of frequency stabilisation)**. A system $S_C$, or an outcome/event $R_j$ occurring in the system, possesses the property of frequency stabilisation (characterised by precision δ) *iff*

(i) it is possible to repeat n 'identical' experiments (with n a number sufficiently large for the precision δ) on the test object T of $S_C$ by observing or measuring T n times under 'identical' conditions C (materialised by the initiating and probing subsystems, or more generally, by the environment); and

(ii) in this experimental series the relative frequencies of the outcomes $R_j$ (j = 1,…, J) of the experiments on T (detected by the probing subsystem) converge towards constant numbers when n grows; more precisely, when n grows the ratio {(number of times that outcome $R_j$ occurs in the experimental series) / n } converges (within δ) towards a constant number, for all j.

**DEF2 (of probability).** Only for an outcome/event $R_j$, or system $S_C$, that shows frequency stabilisation according to **DEF1**, as can be tested by an experimental series in well-defined conditions C, the probability of the outcome/event $R_j$ ($j = 1,\ldots, J$) is defined (within precision $\delta$, and for all j), and given by:

$P(R_j)$ = (number of times that outcome/event $R_j$ occurs in the series) / n (for n large enough to attain frequency stabilisation within $\delta$).

These are, in principle, definitions of 'frequency stabilisation within the precision $\delta$' and of '(empirical) probability-cum-precision-$\delta$'. In **DEF1**, "a ratio converges within $\delta$ to a constant number (P)" means that the ratio varies only in the interval $[P-\delta, P+\delta]$. I believe that the notion of precision is an intrinsic and unavoidable ingredient of probability as empirical notion, as elaborated in **C4** below. However, if one wishes, this model allows one to define 'the' probability, in an absolute sense: it is the number in **DEF2** corresponding to ($\delta \rightarrow 0$, $n \rightarrow \infty$), the theoretical number that can *sometimes* be predicted exactly for idealised systems, e.g. chance games and quantum systems. Hence, one can recover von Mises' limiting frequency definition as a limit case of ours.

Regarding terminology, I will say that the outcome $R_j$ 'has' the probability $P(R_j)$ with or within the precision $\delta$. Note that the outcomes $\{R_j\}$ are the values of a property (variable) X of T probed in the experiment C. Therefore probability can also be associated with (physical) properties/variables (as done in physics, and as can trivially be done for any probabilistic system[17]). Note also that these conventions encompass the usual habit of associating probabilities with 'events' ($E_j$): the event $E_j$ consists in the variable X assuming a value $R_j$ in the experiment C; hence, $P(E_j) = P(R_j)$. **DEF2** implies that probabilities only make sense, are only defined, for systems $S_C$ that exhibit frequency stabilisation in some property X according to **DEF1**. One can also say that the experimental outcome $R_j$ (or even the object T in conditions C) exhibits frequency stabilisation.

I conjecture that the above definitions cover all events or systems that are accurately, numerically described by probability calculus. My main argument lies in the unifying and therefore problem-solving capacity of the model based on **DEF1&2**, as exposed in the next section (results **R1 – R4**) and in auxiliary comments **C1 – C6** below. I cannot review here all probabilistic systems and show how to apply the definitions to them, but I will give examples for several more (classes of) systems, besides the ones already treated; and one can find many more in the works of von Mises (e.g. 1928/1981, 1964).

[17] For instance, for usual die throws the property of the die that is tested/observed is X = number of eyes on the upper face (at rest).

**C1. Compliance with Kolmogorov's axioms**. It is well-known that the frequencies or ratios defining probability in a model like the above obey the axioms of Kolmogorov's probability theory (e.g. von Mises (1964), van Fraassen (1980), Gillies (2000), Khrennikov (2008)). Thus, in von Mises-style theories probability is defined in terms of a measurable quantity; and the axioms of probability theory become theorems in his theory. Most other interpretations of probability do *not* explain why probabilities satisfy these axioms – so this is a genuine argument in favour of von Mises-style interpretations.

**C2. The gist of the interpretation**. **On the ontic (objective) and epistemic (subjective) aspects of probability theory**. If one reduces the definitions to their essence, they state that probability *only* exists for ρ-random systems that *can* be subjected to massively repeated experiments or observations, occurring in scientifically well-defined conditions; and that in that case probability is given by a simple ratio or frequency. Under the scientific realism I assume, probabilities also make sense ('exist') – more precisely: describe something real and objective pertaining to systems – even if no experiments have been or will be done on them; *but it must be possible that such an experiment be done,* at least in principle, as stipulated in **DEF1.** This is the gist of the interpretation. The 'in principle' proviso means: empirical verification of the probability value must be possible, but it can be an *indirect* verification (via a scientific theory): see the discussion below.

Hence, according to this interpretation probability is not a property of individual objects or events *per se*; rather, it is a property of *repeatable experiments* (C above), or of repeatable events, even better of *(collections of) repeatable systems-in-context (ρ-systems)* – a phrasing that resonates with realism. As a side-remark: within scientific realism, the distinction 'object versus composed object (system)' is not utterly fundamental, once one is aware of it.

In the case of physics, the idea that probabilities refer to experiments has been highlighted by several philosophers, including Popper[18] (1957, p. 67) and van Fraassen (1980, Ch. 6); it is implicit in von Mises' theory and Kolmogorov explicitly mentions the experimental conditions, as recalled above (1933/1956, p. 3). I believe this conclusion, that probability is a *property of systems-in-context that can repeatedly be tested in experiments*, should be generalised to any scientific notion of probability, also beyond physics. The general wisdom here is clear: untestable notions are unscientific.

Notwithstanding the operationalist concepts used in the definitions, under scientific realism any sound scientific model and definition should say something about the objective world; and indeed probability values can objectively (or inter-subjectively) be verified and agreed-upon, as all scientific properties can. So each probability (value) must correspond to something real, ontic; if one wishes to give a new name to this objective property, only to emphasise that it is an ontic category, one could call it 'objective chance'. The key idea is: *probability values also govern probabilistic systems when nobody looks*, e.g. molecules in a gas, biological populations, quantum systems, etc. The temperature of a gas, for instance, is a verifiable property of a massive collection of molecules. This temperature is determined by the probability distribution of the molecules' velocity magnitudes, a probability distribution that is an objective feature characterising the movement of the molecules (cf. von Mises (1928/1981 p. 20) and Appendix 4.2). It makes sense to call this feature the objective chance (some prefer 'propensity') to assume given velocity values in given conditions. A sure witness of the objectivity of such a probability distribution is the fact that everyone can measure the same temperature of the gas; neither temperature nor probabilities are subjective on this very science-compatible interpretation.

Indeed, a merit of the definitions is that they apply both to artificial systems (chance games) and natural ones: the latter behave in a ρ-random manner without human intervention. In the case of chance games the initialising and probing subsystems can be identified as truly separate systems acted upon, or made, by humans (more in **C6**: Laplace's interpretation for chance games resorts under **DEF1&2**). In the case of natural random phenomena the initialising and probing systems coincide with – are – the environment, which materialises the experimental conditions mentioned in the definitions. Natural random events are normally not initiated nor measured by humans in repeated experiments, but according to the above model they only have a well-defined probability if they occur, or can occur, on a massive scale as could be imitated in the lab, or studied 'out there' under reproducible conditions. Again, the probability values of such natural events can only be revealed, in any case 'proven', by experiments in their environment or imitating their environment; they can sometimes be hypothesised by calculations based on experimental data or theoretical assumptions[19] (for some examples from physics, cf. Appendix 4.2).

---

[18] Popper (1957, p. 67) famously interprets probability as the *propensity* that a certain experimental set-up has to generate certain frequencies.

Elaborating on this last remark on 'theoretical' probabilities: in view of their primordial link with experiment, some might worry that **DEF1&2** are not applicable to certain probabilities that are predicted by theories and seemingly not directly measurable by experiments. An example is again the probability of the velocity of gas molecules. Since von Mises has treated this case (1928/1981 p. 20), I will only briefly summarise it in Appendix 4.2. The key idea is this: even if experiments for measuring molecular velocities are not available (yet), they are a conceptually immediate extrapolation of possible experiments (they are valid thought experiments), and, most importantly, *the inferences of such thought experiments can be tested by experiments in agreement with **DEF1&2***, as explained in the Appendix. In sum, when I stipulated above that "it must in principle be possible that a repeated experiment be done", the 'in principle' means that the assumption that the property in question is characterised by a probability must be an assumption within a theory that, as a whole, has empirically verifiable consequences. (We know at least since Quine that statements, e.g. attributing probability values, are always tested as part of theories.)

Also, note that the 'repeatable experiment' in **DEF1&2** can be rudimentary: it can consist of simply repeatedly 'observing' a property on an ensemble of test objects 'out there', in their natural environment – a simple type of measurement, simpler than for instance quantum measurements done in sophisticated lab conditions. This possibility is explicitly added in **DEF1&2** to include probabilities extracted from statistical data, describing for instance properties of biological populations, meteorological events, or cosmological systems. Such data can represent the outcomes of tests, observations, done in homogeneous-enough conditions on the elements of the population. From such experimental data probabilities can be extracted – *but only if the data are characterised by frequency stabilisation within a certain acceptable precision*. So, for instance cosmological probabilities studied by Beisbart (Beisbart and Hartmann 2011, Ch. 6) are covered by **DEF1&2**, thanks to the notion of precision (cf. Appendix 4.2 and **C4** below).

---

[19] Von Mises found it important to emphasize that when one applies probability theory to a concrete problem, one always starts from assuming a certain probability distribution. One of his well-known maxims is "Probability in, probability out". E.g., when one calculates the chances in urn picking by combinatorics, one typically (and often implicitly) starts from assuming that the chance that a given ball is picked is equal for all balls. Then one proceeds to the problem posed (such as: what is the probability that two successively picked balls are white ?). In theoretical physics too, one often guesses or has evidence or derives by using physical theory that the probability distribution of a certain property is of a certain type (say a Gaussian). Based on this hypothesis one then uses probability calculus or statistics to derive other probabilistic predictions – where it is often only the latter that can be tested. This is an example of *indirect* verification of probabilities.

A recurrent question that haunts this chapter is to what extent the operationalist concepts in the definitions point to subjectivism or operationalism – *are these fundamental*? It is very intriguing that the model overtly posits an intimate link between probability measures and experiments exhibiting ρ-randomness, or if one prefers between probability and experimental conditions (C). I believe one cannot avoid this operationalist or anthropocentric move; and one wonders whether it has essential philosophical implications, pointing to the idea that (the application of) probability theory has an epistemic, conventional aspect to it. The notion of precision in **DEF1&2** seems to corroborate this idea. I have already argued in Section 3 that such an operationalist / epistemic / subjective element should not frighten us too much; it is ubiquitous in physics.

Indeed, consider again the idea that the probability of X can only be well-defined if well-defined experiments for X can be devised – described by a protocol usable by trained experimenters. This claim seems to suggest that the concept of probability points to our human cognitive limitations. The phrasing in clause (i) of **DEF1** "it is possible to repeat n identical experiments" might trigger a naïve question: *possible for who?* A more skilled experimenter might perform experiments and write protocols disclosing probabilities that remain unknown to others – recall e.g. the elusive hydrodynamic probabilities mentioned in Section 3. Do probabilities only exist for these experienced experimentalists? In a trivial sense, yes, but in a more profound sense, no, of course. "Possible" means "possible for anyone who can", "possible for the most capable experimenter", and ultimately "possible in the absolute". We get here very close to Laplace's famous conjecture that we only need probability because of our cognitive limitations. An all-knowing being would know the course of events with certainty – attributing only zeros and ones to probabilities, as is fitting for a fully deterministic universe. While *we* can only fix and know the (initiating) conditions of a die-throwing experiment up to a certain point, allowing us to extract non-trivial probabilities, one can perhaps legitimately imagine an all-knowing being who knows and fixes the experimental conditions of another, more fine-grained die-throwing experiment, leading to trivial probabilities. Yes, it seems so, *but in any case that would be another experiment, characterised by different experimental conditions*. Which probabilities are correct ? Both – *both refer to different ρ-systems*. In sum, whether one considers the ultimate nature of an event/system as genuinely random as opposed to deterministic, depends on one's knowledge state – on whether one factually knows (or assumes the existence of) an experiment that leads to trivial probabilities (0 or 1) or not (recall that deterministic and genuinely random systems are two sides of the same coin: both are ρ-random). *I would consider this as the most relevant epistemic, subjective dimension of* the application of *probability theory*. An example of how this relativity of (in)deterministic ascription manifests itself is given in Appendix 4.3.

Thus, from the point of view of the philosophy of physics, the mentioned link probability – experiment and the said conventional dimension lead to following fascinating metaphysical question: Are indeterministic (i.e. non-trivially probabilistic) events *ultimately* non-trivially ρ-random, or is there *ultimately* a hidden determinism in nature? Do indeterministic properties, including quantum ones, all emerge from deterministic properties – as Laplace conjectured? This question goes much beyond the scope of this chapter, but, of course, the orthodox interpretation of quantum mechanics claims that indeterminism is ontologically fundamental, and hidden determinism a fantasy. However, I believe this question of the ultimate nature of indeterministic phenomena is an unsolved problem. Against the orthodox position in the quantum physics community, there are cogent and rarely discussed arguments in favour of hidden determinism (Wuethrich 2011, Vervoort 2019, 2020, Nikolaev and Vervoort 2023). I repeat that this is not the same claim as the main contention of the subjective interpretation of probability, which makes probability an epistemic, subjective category. Probability values are objective (cf. **R3** below); but at the same time the question of whether a system is considered indeterministic (non-trivially probabilistic) or deterministic does depend on the knowledge state of the agent attributing probability values. So it is this ascription of a deterministic or indeterministic nature to a system that could be termed epistemic[20] (an idea illustrated in Appendix 4.3).

[20] One easily gets in a Babylonian confusion of tongues in this debate, so I will be explicit, at the risk of repeating myself. In contrast to the general literature on the philosophy of probability, in the philosophy of physics literature the notion of 'subjective' probability is often associated with determinism, and 'objective' probability with indeterminism (e.g. Wuethrich 2011). The underlying thought is that 'objective' probability is the only possibility we are left with in a fundamentally indeterministic (say quantum) world; while conceiving probability 'subjectively', as an expression of our ignorance of underlying deterministic processes, is appropriate in a deterministic world. (In this problem of the ultimate (in)deterministic nature of elementary particles, instead of 'objective' probability I prefer 'ontic'; and instead of 'subjective' I prefer 'epistemic'.) Now, according to my model, even in the latter subjective or rather deterministic

Now, this epistemic claim can, in fact, easily be grounded on an objective basis; it appears to be not fundamental. As already suggested above, it suffices to realise that probability does not belong to systems per se, but to ρ-systems. *It may be undecided whether 'a system' is deterministic or indeterministic, but whether a ρ-system is deterministic or not is an entirely objective matter*[21]. Interestingly, as the example in the appendix illustrates, one always comes to this conclusion: the more knowledge observers have about a system, the surer they will consider it deterministic. The all-knowing observer is a determinist, so to speak. Does this hint to the heretic conclusion that Laplace and Einstein win against Bohr (Vervoort 2019, 2020, Nikolaev and Vervoort 2023)?

But let us not get carried away. These matters of the philosophy and foundations of quantum physics, however interesting, should not divert us from our primary goal, which is to define probability. Further conclusions about the subjective-objective debate are given in the next section (**R3** on the subjective interpretation of probability) and in Section 6, the Conclusion.

**C3. On 'identical' or 'similar enough' conditions**. **DEF1** of ρ-randomness relies on the notion of 'identical' or 'similar-enough' conditions and objects. It may thus look suspicious in philosophers' eyes. But it is not difficult to see that the definition is sound in the scientific sense: it allows to define ρ-randomness and probability in a verifiable manner, a manner that allows for objective or if one prefers inter-subjective agreement. Indeed, the conditions for doing the 'frequency stabilisation test' can be put down in a protocol and communicated for any particular system: "do such-and-such ('identical') initial and final actions on such-and-such ('identical') systems (or put or observe these 'identical' systems in these 'identical' environmental conditions) – and the probabilities $P_j$ will emerge. I found frequency stabilisation and the probabilities $P_j$, so if you do the experiment in the 'right, identical' conditions, you should find them too." The term 'identical / similar', which is ubiquitous in the philosophy of probability theory, can thus be

---

or epistemic view of probability, probability values are still objective measures in the strict sense (as used in the philosophy of probability) – namely relative frequencies. So it is more precise to consider not the probability value itself as subjective, but rather the "ascription of a deterministic or indeterministic character" to the event under scrutiny (cf. Appendix 4.3, based on Vervoort 2020). For a related view, see (Bacciagaluppi 2020).

[21] Similarly, the question of whether a variable/property X of system T in conditions C is deterministic or indeterministic has an objective answer, in principle: it depends on whether a repeatable experiment (measuring X and representing C) on T exists that always leads to trivial probabilities for the values of X. If so, X (and T in conditions C) are deterministic. For instance, the spin $\sigma_Z$ of an electron T in a state $\psi = (u_z + u_{-z})/\sqrt{2}$ measured along the axis z is *indeterministic* (the context C is here described by the electron's initial state ψ): the spin has two possible values (+1 and -1) and both have a probability ½.

defined in an operationally consistent manner: *'identical' is 'identical' insofar a probability distribution ρ emerges*.

Thus, it is helpful to introduce the notion of 'p-identical', or rather 'ρ-identical', events / systems. To identify the probability of an event (R), one needs to perform a series of experiments on ρ-identical systems (including test object, initiating and probing subsystem, or environment) – systems that lead to the probability distribution ρ(R) of the results R generated by the system. Vitally, all subsystems *combined* lead to a stable ρ; ρ is only well-defined with respect to such a partitioned system in which all parts 'work together' to lead to precisely ρ. Changing one element of the system will lead, in general, to a different ρ. Therefore I define the concept of ρ-system = a system {test object, environment, initiating/probing subsystem} that is characterised by frequency stabilisation (in at least one of its properties). Many problems can be avoided by using the notion of ρ-system, or ρ-identical system, for instance the 'reference class problem' treated in **C4** below and other problems addressed in Section 5. Again, in a slogan: probability belongs to ρ-systems repeated/repeatable in ρ-identical conditions.

**C4. On the meaning of 'n' and the role of 'δ'. A difference with von Mises?** Frequency stabilisation and probability are defined by means of convergence of a certain ratio when n grows. Here n is a variable or index that counts the number of trials or repetitions of an experiment; more precisely, it represents the cardinality of a set of outcomes of a real experiment that are subsequent in time (or of a set of outcomes randomly drawn from the experimental series, cf. von Mises 1928/1981, 1964). As most experimentally determined quantities, an experimental probability measure has a precision δ; and scientific practice, especially in physics, has provided overwhelming evidence that this precision δ is a function of n, as noted above. The link δ ↔ n is ultimately a (hypothesised) fact of nature, and an entrenched tool of science and statistics.

Clearly, the phrasing in **DEF1&2** is close to von Mises' original definition of probability as $P(R_j) = \lim_{n\to\infty} n(R_j)/n$. So close that one may wonder whether there is a real difference. Sure, the definition "P(Rj) = n(Rj)/n for a number of trials n that is large enough for the precision δ" avoids the notion of infinity which some find confusing (recall Appendix 4.1 and Section 2). But note that from a pragmatic point of view, if one would use **DEF2** to experimentally determine a probability *and if the ρ-system is repeatable 'at will'*, my definition leads to numbers that are equal, *to any desired precision,* to those identified via von Mises' definition. At least operationally there is then no difference in the definitions: they can lead to the same results, the same numbers. Scientifically speaking they are identical, one might say.

But at least on the surface my definition is more inclusive: it allows to determine 'probabilities with a precision $\delta$', or $\delta$-probabilities, even when one cannot determine n(Rj)/n for n so large that $\delta \to 0$, which happens in many somewhat complex systems, e.g. from meteorology, biology, etc. For 'simple' physical systems as a die or as certain quantum systems, one knows the conditions for repeating the series (quasi-)perfectly; one knows the $\rho$-system (quasi-)perfectly; and one can determine probabilities with practically $n \to \infty$ and $\delta \to 0$. For other systems one has no access to massive data, for instance because one cannot repeat the conditions C (the $\rho$-system) very many times.

However, I claim that even in that case von Mises' definition as limit is entirely sound as a construct of a natural science theory: it refers to the mathematical object one *would* identify if one *could* repeat C infinitely many times. I have already argued in Section 3 that such an expression is not different-in-kind from other counterfactuals met in science; in natural science mathematical objects, say a field-value from electrodynamics, are always compared to their experimental real-world referents under counterfactually idealising conditions ("this mathematical object, this number, would be found/measured if these idealised, counterfactual conditions would be met"). In other words, mathematics-based theories in physics refer to quantities that are always only approximately measurable, quantities of which the real-world referents always depend on experimental conditions C, and that are (almost) always characterised by a precision $\delta$. So one should not expect that a real-world (experimental) probability, the object of von Mises' theory, *exactly* corresponds to a theoretical probability, as they can *in rare cases* be calculated with infinite precision ($\delta = 0$) within certain assumptions, say by Laplace's formula or by quantum mechanics. Demanding that von Mises' probability should be verifiable with infinite precision is as unreasonable as asking that any physical property be measurable with infinite precision.

The often-heard complaint "frequencies do not always stabilise" reveals a confusion of the same kind: only in idealised cases can one repeat *exactly identical* conditions C at will; often these conditions (describable by a protocol) are only approximately known and controllable – except in simple cases. It is again here that my definition allows to determine approximate probabilities corresponding to frequencies that only approximately converge (within $\delta$) *due to only* *approximately fixed conditions*. Indeed, if some frequency does not converge with $\delta \rightarrow 0$, even on sufficient repetitions, my model suggests that this is due to some insufficiently fixed / stabilised condition-variable(s), i.e. some of the variables in the sets $C_i$, $C_f$, $C_e$ in $P(R|C_i, C_f, C_e)$, with $C_i$, $C_f$, $C_e$ representing respectively the initiating, probing and environmental conditions. (As an example, recall the case of the variation in magnetic field given in Section 3, due to a varying systematic cause in $C_e$.) Indeed, as we will elaborate in the next chapter, these conditional variables are in *causal relationship* with R (in agreement with our intuition that these conditions co-determine R). Therefore, fluctuations in some of these variables should lead to fluctuations in R and hence in P(R) (= $P(R|C_i, C_f, C_e)$). In the last discussion I have assumed that, even if $\delta$ does not go to zero, it is sufficiently small to speak of the probability of the variable at stake (R). This means that the most important causal factors (conditions) are fixed / constant; the fluctuating variables have a much weaker influence on R. I have discussed these issues here without going into the mathematics; I conjecture that such an exercise would be interesting and could, perhaps, lead to quantitative insights.

In conclusion, one can defuse criticisms against the limit-definition if one understands it correctly as referring to the mathematical object one *would* determine if one *could* repeat *exactly* the conditions C infinitely many times. My definition boils down to a more explicit variant.

In his (2009) Hájek considers an experiment of coin-tossing in a train that moves periodically back and forth on its track pointing west-east (see details in (2009), p. 218). The outcomes are recorded at well-defined space-time points, and can therefore be represented in a space-time diagram. In the case considered they form the following periodic series as a function of time: HHTHHTHHT…, so that the probability of obtaining H, P(H), = 2/3. The author claims that one could also 'look at the data' on the space-time diagram from west to east, and then one 'sees' the series HTHTHTH…. with limiting frequency P(H) = ½. The author concludes: "Now, why should one answer have priority over the other? In other words, we have more than one limiting relative frequency, depending on which spatio-temporal dimension we privilege." But this objection does obviously not apply to von Mises' frequentism, which explicitly attributes probabilities to real experimental series, not of course to deterministically selected partitions of outcomes, for example by projecting them on the spatial axis in a space-time diagram! So, 'n' in the definitions **DEF1&2** is a number that represents the cardinality of a series of results that follow each other in time in the experimental series (or of a series of outcomes randomly drawn from the experimental series, cf. von Mises 1928/1981, 1964). [Hájek's example is contrived also because the experimental series he considers is perfectly periodic, i.e. deterministic: there must be a mechanism ensuring such a perfect periodicity. A toss in the forward moving train always gives H, in a backwards move always T; there are two times more forward moves. So the 'probability' 2/3 only makes sense for a fully unaware observer (a robot), so a ρ-system that *only* includes an observing system registering the outcomes, not a ρ-system that includes a movement detector – which would allow one to define and measure probabilities P(H|m) = 1 for m indicating that the die toss is done in a forward move.]

In sum, objections as to an alleged problem of undefined limiting values (Hájek 2009, 2023) are dissolved by realising that sets of real experimental series must be considered. It should be noted that several criticisms against the frequency interpretation neglect the real-world basis of experimental series of frequencies: these are physical or real-world series, not logical nor mathematical ones (as another manifestation of this confusion, see the criticism by Fine (1973), p. 93, and my objection in Appendix 4.1). The 'reference class problem' (Hájek 2009, 2023) is solved immediately in a similar way: probabilities are relative to a well-defined experiment, or if one prefers to a ρ-system, or if one prefers to the outcomes – ordered by increasing n – of well-defined experiments involving ρ-systems.

(I hope that the reader, by the end of this chapter, will agree that the operationalist / epistemic elements in probability theory are innocuous, but for some philosophers a definition closer to von Mises', not relying on precision, might still be more satisfactory, *metaphysically speaking*: one could say that 'the real' probability is the limiting frequency 'at infinity' – if one is aware of the surprisingly many implicit notions discussed above. An ideal, theoretical number, corresponding to the theoretical precision interval $\delta = 0$, a number that one can rarely know with infinite precision but that one can determine with increasing precision when n grows in stable conditions.)

**C5. Applicability to continuous variables**. In physics most ρ-random phenomena are characterised by continuous variables rather than discrete ones. Kolmogorov was the first to rigorously show that his calculus applies to continuous variables too (Kolmogorov 1933/1956). Thus, in this case the variable R will range over an interval of the real numbers. A notorious exception are certain quantum systems, notably when they are carefully 'prepared' – I said

'initiated' above. Still, many or most quantum properties have a continuous spectrum. Note that the origin of the discreteness of outcomes of chance games is obvious: the highly symmetrical objects used are constructed on purpose by humans to generate a limited number of outcomes.

In the case of a continuous ρ-random variable R one defines in probability calculus the probability density function ρ(R), the meaning of which is however defined via the concept of probability. Indeed, in the continuous case the probability P(R is in dR around the value r) = ρ(R = r).dR. This equality defines ρ via P as characterized in **DEF2** (here P(R is in dR around r) is the probability that the variable R assumes a value in the infinitesimal interval dR centred on R = r). One can thus safely treat probability density functions and probability on the same foot; a model that interprets probability for discrete variables also does the job for continuous variables. Note that one can formally describe discrete cases by a density function $\rho(R) = \Sigma_j P_j.\delta(R\text{-}R_j)$, where δ(−) is the Dirac delta-function: the discrete case is a special case of the more general continuous case.

**C6. Unification with the classical interpretation**. Frequentist models in general cover the classical interpretation of probability, due to Laplace and others, which is used to tackle chance games. Since von Mises (1928/1981, p. 66 ff.) and many modern texts have proven this result, I will briefly present an argument as applied to **DEF1&2**. Notice, first, that **DEF1&2** *can* be applied to such chance games. The initiating / initialising subsystem is most of the time a randomizing hand (tossing a coin or a die, pulling a card or a ball from an urn, etc.); the probing subsystem is often simply a table (observed by a human).

In Laplace's well-known terminology, for calculating a probability of a certain outcome, one should consider "events of the same kind" one is "equally undecided about"[22]; within this set, the probability of an outcome is the ratio of "favourable cases to all possible cases". Notice now that dice, card decks, urns containing balls, roulette wheels, etc. are constructed so that they can be used to produce *equiprobable* (and mutually exclusive and discrete) *basic outcomes*, i.e. *having all $P_j$ equal, and given by 1/J* (J = the number of basic[23] outcomes). Equiprobability is at any rate the assumption one starts from for making mathematical predictions, and for playing and betting (recall that one has to make *some* probabilistic assumption to start from, according to von Mises' maxim “probability in, probability out”); and indeed Laplace’s “events of the same kind one is equally undecided about” would now be termed equiprobable events. Along the lines exposed above, a chance game can thus be seen to correspond to a ρ-system with $\rho(R) = \Sigma_j (1/J)\, \delta(R-R_j)$.

[22] The events thus fulfil the 'principle of indifference' introduced by Keynes.
[23] The 'basic events' or 'basic outcomes' of coin tossing are: {heads, tails}, of die throwing: {0, 1, …, 6}, etc. The probability of 'non-basic', or rather composed, events (an even number of eyes; a six or an ace; etc.) can be calculated by using probability calculus and combinatorics.

In the special case of equiprobable and mutually exclusive events, the frequentist and classical interpretation indeed lead to the same numerical values of probability. Within the frequentist model, $P(R_j) \longrightarrow \frac{n.\frac{1}{J}}{n} = 1/J$ (the numerator n / J = the number of $R_j$-events among n (>>J) exclusive and equiprobable events each having a probability 1/J). The result, 1/J, is equal to the prediction given by Laplace’s formula (1 favourable case over J possible cases).

Note, however, that Laplace’s formulation is not superfluous. It allows in the special case of chance games for calculation, i.e. theoretical prediction: in complex examples ‘favourable cases’ and ‘possible cases’ can conveniently be calculated by the mathematical branch of combinatorics – a theory of counting, initiated by the fathers of probability theory. In sum, it appears that the classical interpretation is only applicable to a small subset of all probabilistic systems, namely chance games, in general artefacts having high degrees of symmetry leading to easily identifiable, discrete basic outcomes. For this subset Laplace’s interpretation can be seen as a *formula*, handy for calculation, rather than an interpretation of probability.

These notes conclude the basic description of my interpretation of probability. As always, the way to validate and strengthen a model is to show that it applies to non-controversial cases, and that it solves problems left open by other accounts. Besides in the examples given, it is useful to explicitly verify the application of the definitions in the somewhat subtler case of probabilities hypothesised in physics (Beisbart and Hartmann 2011, pp. 143 – 167). A few more examples are given in Appendix 4.2; the examples given throughout this chapter hopefully illustrate how to apply the definitions also to other cases. The exercise in the appendix illustrates, notably, the essential role of the ‘environment’ (say temperature), part of the ρ-system, in the case of natural random phenomena.

## 5. More problems solved

In this section I show that the account based on **DEF1&2** allows us to interpret controversies and to solve a variety of problems of the interpretation of probability. Some conundrums have already been addressed in the previous section. Here I will propose my solutions to problems **R1** – **R4**, on single-case probabilities and the unification with the propensity interpretation (**R1**); on Bertrand's paradox (**R2**); on the unification with certain subjective interpretations and the origin of the subjective shift (**R3**); and on the probabilistic origin of many quantum puzzles (**R4**).

**R1**. **On single-case probabilities and the unification with the propensity interpretation**. An immediate consequence of the model based on **DEF1&2** is that it dissolves the intensely debated problem of 'single-case probabilities' (see references in Beisbart and Hartmann 2011, Ch. 1 and a detailed discussion in Gillies 2000, Ch. 6-7). According to von Mises-type interpretations the solution is clear: they imply that it makes no sense to talk about the probability (as a quantitative measure) of an event that cannot be repeated in experiments. For instance, it makes no sense to speak about the 'probability' of a dictator starting a world war and the like: no experiments can be repeated here, and even less experiments in well-defined conditions. One might of course attribute a subjective *likelihood* to such events or hypotheses; but we are interested here in the scientific notion of probability. No experiments, no natural science. (As already hinted to above, one might perhaps still speculatively use the notion of probability *in philosophy* to refer to hypothetical properties of classes of dictators that would share enough similarities to envisage stable statistics. Perhaps this speculation has some use in certain philosophical contexts, cf. § 6.2, but surely no quantitative information could be extracted from it.)

Popper's propensity interpretation of probability was an attempt to accommodate such probabilities of single events – quantum events were his inspiration, such as the disintegration of one particular atom (Gillies 2000, Ch. 6). Intriguingly, previously Popper had done extensive work on the frequency interpretation, but he ultimately rejected the latter in order to make sense of single-event probabilities (Popper 1959). Now, this inspiration from quantum mechanics is quite puzzling: what is the difference with classical systems? What is the relevant difference between the probability of this particular die showing a 6, and the probability of this particular atom disintegrating in the next hour? I do not know of any relevant difference. The practice of physics has shown that *any probabilistic quantum property of any single quantum system* can only be attributed a probability value if the property can be measured on an ensemble of such systems – as in the case of macroscopic systems. Sure, one can *metaphorically* say that *this* atom has a certain disintegration probability, just as one can metaphorically say that this die has a certain probability to show a six. But this has no immediate empirical meaning other than the meaning conferred by the frequentist account. And indeed, propensity interpretations typically state something along these lines: a probability P of an outcome/event is the propensity of conditions C, or of repeatable experiments, to produce outcomes/events characterised by a frequency P. Such a definition is still based on a frequency interpretation. Therefore, the added value of the propensity interpretation seems far from obvious. Popper has argued that frequentism associated with conditions, as they are an explicit part of **DEF1&2**, amounts to a propensity interpretation. But since a well-understood frequentism answers so many more questions, as argued in this chapter, this is turning the logic upside-down: it is the propensity interpretation that could be seen as a trimmed-down variant of a detailed frequency interpretation.

Perhaps proponents of the propensity interpretation prefer it because the expression 'probability is a propensity of conditions' has a more evident objective touch to it; it confers the impression that there is *something in nature* that creates probabilities. Popper compares propensities to physical forces (Popper 1959). That is fine in itself; but I fear we enter here in the realm of semantic fuzziness and that there is nothing substantial to be gained, that there are no new important questions answered. I argued above in some detail that my frequency interpretation also comes to the conclusion that probability is an objective category (notably for the obvious reason that frequencies emergent from ρ-systems are objective); it is in any case a basic assumption of my model; and see **R3** below. Also, it is no more objective to attribute probabilities to single atoms than to ρ-systems.

**R2. Solving Bertrand's paradox**. A classic paradox of probability theory, debated by the fathers of probability calculus and many mathematicians since, is Bertrand's paradox. It goes as follows: "A chord is drawn randomly in a circle. What is the probability that it is shorter than the side of the inscribed equilateral triangle ?" (see e.g. Marinoff 1994). Bertrand showed in 1888 that apparently three valid answers can be given – so which is the right one?

Bertrand's problem is an archetypical case of a paradox that originates from neglecting the fact that probability is physical, not only mathematical. Bertrand's puzzle refers to a real-world situation; it is not only a mathematical problem. In the real world there are many ways to 'randomly draw a chord' (which may not be obvious upon first reading of the problem – hence the confusion). One can for instance randomly chose two points (homogeneously distributed) on the circle by using a spinner[24]; a procedure that leads to the probability 1/3, as can be measured and calculated. But two other (classes of) randomisation – *or initiating* – procedures are possible, leading to two different outcomes (1/2 and 1/4). In sum, *Bertrand's problem is not well posed*, a conclusion that is now quite widely accepted (see e.g. Marinoff (1994))[25]. Note that this conclusion immediately follows from **DEF1&2**, according to which probability is only defined for experiments in well-defined conditions, *including initiating conditions* – it is these that are not defined in the phrasing of the paradox. Probabilities need to be associated with a ρ-system, which includes, in the case of chance games, an initiating subsystem. It is not surprising that in a situation that is experimentally ambiguous no unique probability values can be defined.

It will be seen that the literature abounds with paradoxes which stem from formulations that are experimentally complex or ambiguous. The infamous 'Monty Hall' problem, popularised by Marilyn vos Savant, will be seen to become equally straightforward if one remembers that, here to, one looks for a relative frequency of an outcome R of a well-defined experiment[26].

**R3. Unification with certain subjective interpretations. Origin of the subjective interpretation**. The subjective or Bayesian interpretation of probability, due to authors as De Finetti, Ramsey, Jeffreys, Jaynes and several others, comes in as many variants as there are authors (see excellent reviews in von Plato (1994), Gillies (2000)). These range from full-fledged subjective versions such as De Finetti's, interpreting probability as a measure of strictly individual belief, to variants that seem easier to adopt in scientific contexts, such as Jaynes' theory. According to subjective interpretations, depending on the variant, probabilities can be attributed to events, hypotheses, propositions, etc. Neutral assessments often come to the conclusion that the hard-boiled subjectivist versions 'become empty' after detailed scrutiny (Gillies 2000, p. 84).

---

[24] More precisely, by fixing a spinner at the centre of the circle; a pointer on the spinner and two independent spins generate two such independent random points, through which then a 'random chord' can be drawn.
[25] The latter reference also reviews von Mises' and Gnedenko's treatment of the problem, cohering with the one here presented (p. 23). Keynes, for unclear reasons, holds fast to the principle of indifference (p. 23).
[26] R = 'I win a car by switching'; P(R) = 2/3.

In the case of Bayesianism à la Jaynes, the situation is *a priori* less clear. How does such a theory compare to frequentism? This is an important question, because Bayesianism has become popular in several scientific fields, including in quantum foundations and quantum information theory (e.g. Caves et al. 2002, 2007, and Bub 2007). It has already inspired valid scientific results. So, if frequentism is the mother of all interpretations, following question is pressing: can these results also be obtained within the frequency interpretation? I conjecture that they can. Although I cannot review here this broad Bayesian research program, I will show that some variants of Bayesianism can be absorbed in the frequency interpretation. I will focus on Jaynes' theory (1989, 2003), to which e.g. Caves et al. (2002), (2007) and Smerlak and Rovelli (2007) refer. But first I will start by opposing the full-fledged subjectivist interpretation to frequentism; this will suggest a hypothesis of why the 'subjective shift' is so tempting.

It is a popular and very tempting idea that, somehow, "probability depends on our knowledge, or on the subject" (cf. **C2** above). When I throw a regular die, I consider the probability for any particular throw to show a six to be 1/6. But what about my friend Alice who is equipped with a sophisticated camera allowing her to capture an image of the die just before it comes to a halt ? For her the probability seems to be 0 or 1. Or imagine Bob, waiting for a bus, only knowing that there is one bus passing per hour. He might think that the probability he will catch a bus in the next five minutes is low (say 5/60). Alice, who sits on a tower having a look-out over the whole city, might have a much better idea of the probability in case (she might quickly calculate, based on her observations, that it is close to 1). Are these not patent examples of the idea that a same event can be given different probabilities, depending on the knowledge of the subject ? And is in that case probability not a measure of the strength of belief of the subject who attributes the probability values ?

Examples as these are unlimited, but the above claims can both be explained and countered by referring to the boundary conditions that are explicitly part of detailed frequentist interpretations, in this case the conditions of observation. Doing a normal die throw, and observing the result on a table as is usually done, corresponds to a well-defined ρ-system, with a specific initiating system, probing system etc. In the example, observer Alice does not measure the same probability: she does not measure the probability of finding a six on a table after regular throwing and regular observing, but of finding a six after measurement of whether a six will land or not on the table. The latter is a very different, and indeed fully deterministic, experiment; *at any rate, the observing/probing subsystem (including a high-speed camera) is very different*. A similar remark holds for the bus-case; the measurement system (even if just the human eye in a given location) is part of the ρ-system; one cannot equate observer (probing sub-system) Alice and observer Bob if their means of observation are widely different. Ergo, Alice and Bob do not identify probabilities of the same event or system: *that* is why they attribute different numbers.

For readers who are not convinced, let us go over the examples more carefully. One may consider the probability of event $e_1$ = "the outcome R = 6 in a regular die throw (probed as usual on a table)". One could also consider the probability of $e_2$ = "R = 6 after a camera has registered that R = 6 one microsecond before the die came to a halt". The probability of $e_1 = P(e_1) = 1/6$, while $P(e_2) = 1$ (suppose so). Within my model, these probabilities are *not* different because of different subjective knowledge states. First note that also for $e_2$ one can define the probability without referring to 'subjective knowledge' or 'information'; $P(e_2)$ could be measured by an automat[27]. So $P(e_1) \neq P(e_2)$, not because someone has a different strength of belief, but because $e_1 \neq e_2$: both probabilities concern different events, different experiments, and in particular different observing conditions and systems. Every probability can be seen as the result of a series of automated tests.

Now, since nothing is simple in the foundations of probability, it is at the same time true that $P(e_2)$ can also be expressed in a manner that invokes, implicitly or explicitly, the status of knowledge of some observer or 'probability attributor'. $P(e_2)$ can be considered a conditional probability, namely the probability that R = 6 *if it is given* that R = 6 one μsec before halting; in symbols $P(e_2) = P(e_1 \mid R = 6$ one μsec before halting$)$. *And the phrasing "if it is given that" seems indeed equivalent, in this context, to "if it is known that"*. This semantic equivalence may well explain why one is tempted to interpret conditional probability as dependent on the knowledge of the probability attributor[28]. On the view I defend here, such a subjective shift could be innocuous as long as one realises that it is a shortcut of thought – it is by no means the only interpretation, let alone the most encompassing one. Therefore, based on the explanatory power of **DEF1&2**, the natural inference is that subjective interpretations, attributing a probability value P(p) to a proposition 'p', only have a coherent content if there is an unambiguous map from P(p) to P(e) as defined in **DEF1&2**, with 'e' an event occurring in a repeatable experiment (or if P(e) is indirectly testable). Needless to say, this condition greatly restricts the use of 'probability' as a scientific notion. The idea that it could be generalised as a philosophical, speculative notion is briefly entertained in Section 6.

[27] Such an automated experiment is this: let a robot launch a die, let it select by camera vision those trials that show R = 6 one μsec before the die comes to a halt, and let it measure on that ensemble R again at full stop of the die (all this could be done by a machine). The relative frequency of these results will converge to 1, as our robot could determine.
[28] Van Fraassen (1980, p. 164ff) gives a detailed analysis of how subjective interpretations (probability linked to ignorance) originating from chance games slip into the objective statistics of physics. Also Gnedenko gives an enlightening analysis of this subjective shift (Gnedenko (1968), p. 26ff).

It is important to emphasise that conditional probabilities can be interpreted in an objective way, since they play an essential role in probability theory; on my view they are the basic notion. Recall that I stipulated that it is often useful to replace P(R) by P(R | $C_i$, $C_f$, $C_e$) where $C_i$, $C_f$, $C_e$ are all relevant variables that describe the initial, final, and environmental experimental conditions, respectively. (As noted above, this replacement can only be done if the C-variables are constant for all trials in the experiment at stake: they must define this experiment.) Von Mises emphasised, and I agree, that *all* conditional probabilities can be seen as corresponding to series of automated experiments in which no intervention nor belief state of a human agent is needed, whether these probabilities are related to chance games as in the above case, or to natural phenomena (1964, pp. 22-24). But is this not a happy argument for the homogeneity of probabilistic phenomena? Under the usual interpretation all natural stochastic phenomena occur according to probabilistic laws also without a human pondering about them. Therefore, in accordance with an acclaimed rule of epistemic parsimony, this gives a head start to interpretations in which the same holds for human-induced events such as die throws.

To corroborate, let us have a closer look at Bayesianism as taught by Jaynes (1989, 2003). In Jaynes' words (2003, p. 88): "We want probability theory to indicate which of a given set of hypotheses {$H_1$,$H_2$,…} is most likely to be true in the light of the data and any other evidence at hand". In Bayesianism the key role is played by Bayes' formula, which Jaynes writes as follows (2003, p. 89):

$$P(H|DX) = P(H|X).P(D|HX) / P(D|X),$$

with X = prior information, H = some hypothesis to be tested, D = the data. The formula allows to derive the 'posterior' probability of H in the light of new data or information, P(H|DX), as a function of its 'prior' probability, P(H|X). "Equation (4.3) [the above formula] is then the

fundamental principle underlying a wide class of scientific inferences in which we try to draw conclusions from data" (2003, p. 89). Thus, (this type of) Bayesianism understands probabilities as attached to hypotheses of which the likelihood can be calculated in view of 'background information, data or evidence'.

Now, if one investigates the concrete cases to which the above formalism is applied, there is little doubt that these can be recast in the vocabulary of the frequency interpretation. The translation would be immediate if in all cases 'H' and 'D' appearing in the Bayesian probabilities P(H|DX) or P(D|HX) would correspond to *outcomes of a random experiment*, and if X would describe the precise *initial, final and environmental conditions* ('C' in my notation). Then the meaning of the above formula would simply be the same as derived within my von Mises-type model. In all examples I examined, this mapping indeed appeared to be possible. Let me show it in detail for one example, the one that Jaynes treats first and in most detail (2003, p. 93 – 96). In this example we have:

X = 11 machines fabricate widgets, which pour out of each machine into a box (1 box per machine); 10 of these 11 machines produce 1 defective widget out of 6 (these are 'good' machines); one machine (the 'bad' one) produces 1 defective (bad) widget out of 3. We randomly choose one of the 11 boxes without identifying from which machine, and randomly pick a widget from it. Take then:

D = The picked widget is bad.

H1 = We chose a bad box.

H2 = We chose a good box.

Jaynes then calculates probabilities as P(H1|D), the probability that 'hypothesis H1 is true if data D is known'. It is however clear that all the calculated probabilities (p. 93 – 96) can be understood as probabilities of outcomes of an automated experiment, in which a robot randomly picks first a box, then a widget from that box, and finally identifies whether box and widget are good or bad. In the objective phrasing P(H1|D) is the probability of *outcome* H1 (a bad box is picked) given *outcome* D (a bad widget is picked), i.e. the probability that a bad box is picked if it is given that the picked widget is bad, all in the experimental conditions described by X (the result is P(H1|D) = 1/6). Even if Jaynes prefers a subjective vocabulary, this is of course a most classic application of probability theory, which can well be done within a detailed frequency interpretation.

I did not do a systematic study of all Bayesian applications to scientific problems; but in the works I consulted (e.g. examples in Jaynes 1989, 2003, von Plato 1994, Gillies 2000, Tijms 2004), it appears that 'hypotheses' and 'data' or 'evidence' can be mapped to outcomes of experiments; their probabilities can univocally be given a numerical value as a relative frequency in a large set of experimental data. I believe this is not a surprise, nay almost obvious, because, again, scientific concepts need a means of verification via experiments. Therefore I conjecture that my above conclusion can be generalised: *if this link with experiments and **DEF1&2** cannot be made, Bayesianism cannot lead to new results* (e.g. concerning the probability 'that a certain theory is correct'; see Bunge (2008) and Tijms (2004), p. 252 for misuses of Bayesianism).

Occasionally, this is a conclusion to which Jaynes himself tends. When investigating whether Bayesianism can attribute a probability to 'the hypothesis that Newton's theory is right', he concludes: "[…] there is little hope of applying Bayes' theorem to give quantitative results about the relative status of theories" (2003, p. 139). Von Mises would have relinquished hope from the start: the hypothesis 'Newton's theory is right' is not a possible outcome of a random experiment. As a consequence, Bayesianism à la Jaynes, to which recent authors in quantum mechanics refer (e.g. Caves et al. 2002, 2007 and Smerlak and Rovelli 2007), should not be termed subjective, *pace* these authors. In comparison to truly subjective interpretations of probability, it is in fact an objective model, even if phrased in subjective terms. This is again almost explicitly stated by Jaynes himself (2003, p. 45).

In conclusion, I contend that Bayesianism à la Jaynes can be translated into the frequency interpretation: the (background or prior) information/knowledge of the subjective phrasing corresponds to the precise experimental conditions of the frequency interpretation. By extension, I conjecture that all coherent variants of the subjective interpretation can be unified with von Mises-style frequency models. This is still a conjecture, since I did not investigate in detail all these subjective interpretations nor their application to problems. But besides by the reasons given above, this conclusion is corroborated by the fact that frequency models have a much greater explanatory power, as I argue in this chapter. Subjective interpretations face many unanswered questions, and are, by the same token, difficult to unify with natural science. (It may still be that they are a handy short-cut of thought and/or a superior conceptual tool for formalisation for certain *localised* problems; but that is another matter and needs to be investigated in more detail. For only one example of such a localised problem for which the subjective interpretation might be helpful, cf. Uffink 2011.)

**R4. A correct interpretation of probability is essential for understanding key problems of (the interpretation of) quantum physics**[29]. The above conceptual analysis brings in focus a striking similarity between classical and quantum systems. The philosophy of quantum mechanics is a very vast field, so I can here only advance what I believe to be the essential arguments in support of the above claim. I must here necessarily remain succinct.

It belongs to the basic elements of the standard or Copenhagen interpretation of quantum mechanics that "the observer belongs to the quantum system", or at least that the measurement detectors belong to it, and influence it. Let me cite Bohr in his famous debate with Einstein on quantum reality (Bohr 1935). Einstein and two of his collaborators, Podolsky and Rosen (together EPR) had initiated this debate in 1935 (Einstein et al. 1935), with an article that is not only an exceptionally impactful paper in physics (Fine and Ryckman 2020), but also one of the milestones of the philosophy of quantum mechanics. EPR contest in the article that quantum theory is a 'complete' theory: they argue that the probabilistic predictions of quantum mechanics must be explainable, must be completed, by additional variables. In short: quantum probabilities must emerge from a deeper stratum of reality – which is essentially what Laplace thought about probabilities. Niels Bohr, one of the fathers of quantum theory, was not amused, and replied shortly after[30]. The key point of Bohr's counterargument is summarised in following somewhat mystifying line: "The procedure of measurement has an essential influence on the conditions on which the very definition of the physical quantities in question rests" (Bohr 1935, p. 1025). To explain the 'essential influence' that any measurement has, Bohr invokes the quantum of action linked to Plank's constant h: any measurement of any property needs the physical interaction of a detector / observing system with the test object under study, an interaction that carries minimally one quantum of energy, and that potentially exerts a perturbing influence. Bohr's phrase is one of the more famous quotes of the Copenhagen school, but even experts as John Stuart Bell, author of Bell's theorem, have complained that it is incomprehensible (cf. Bell 1981 p. 58, where Bell also identifies the above quote as the key ingredient of Bohr's reply – I agree).

[29] This sub-section is based on older work, where a few more ideas can be found (Vervoort 2011, 2012).

[30] For a detailed interpretation of the EPR article, cf. (Fine and Ryckman 2020). In essence, EPR claim that certain (entangled) quantum states possess 'elements of reality' which cannot be described by quantum theory. Ergo, it must be incomplete (it lacks certain variables describing these elements of reality). In somewhat more detail, and simply said, EPR claim to have demonstrated that non-commuting properties as the position (x) and momentum ($p_x$) of a particle must exist simultaneously – in overt contradiction to Bohr's complementarity principle. Bohr's reply can be summarised thus: x and $p_x$ cannot be measured simultaneously; ergo they do not exist simultaneously; hence quantum theory is complete after all. Who is right? That is still debated, but notice that Bohr's answer is heavily anthropocentric and operationalist/instrumentalist, putting it a priori in serious tension with scientific realism… and therefore possibly with science (recall the MCT-thesis). Let us not say more here.

Incidentally, I could not make sense of Bohr's sentence before I did the present study on probability. My interpretation of probability offers a natural explanation of what Bohr must mean: namely that *(the values of) quantum properties depend in a fundamental way on the observing conditions, or, if one prefers, on (the interaction with) the observing subsystem*. Quantum properties are probabilistic properties; and I argued throughout this chapter that any probability is associated with a ρ-system and depends inextricably on the observing subsystem / conditions. In classical systems one has to look a bit more carefully for examples to exhibit this in-principle fact, but in the quantum realm it becomes of *ubiquitous and manifest importance*. One is tempted to conclude: quantum systems are probabilistic systems *of a pure kind* – more on this in a moment. Just one example: the probability that an x-polarised photon passes a y-polariser obviously depends on x and y; here the angles x and y are the variables that describe the initiating and probing conditions.

A few salient statements of A. Peres' textbook on the foundations of quantum mechanics, nicely expounding the orthodox Copenhagen view, confirm this idea: "A [quantum] state is characterized by the probabilities of the various outcomes of every conceivable test" (Peres 2002, p. 24), where these 'tests' should be understood as pertaining to an experimental set-up. A little further: "Note that the word 'state' does not refer to the photon by itself, but to an entire experimental setup involving macroscopic instruments. This point was emphasized by Bohr […]" (Peres 2002, p. 25).

We are by now ready to conclude that these statements can be generalised to *any* probabilistic system, quantum or classical. To a die throw one could obviously associate a 'state' consisting of the six discrete results and their probabilities. Also the probability of an outcome of a die throw depends crucially on the 'entire experimental setup' that is used to perform it. Compare to the following passage (Peres 2002, p. 73): "The notion of density matrix – just as that of state vector – describes a *preparation procedure*; or, if you prefer, it describes an ensemble of quantum systems, whose statistical properties correspond to the given preparation procedure". I termed the 'preparation' of this statement the 'initiating / initialising' of the probabilistic system. I believe that instrumentalist passages à la Copenhagen of textbooks as (Peres 2002) all can be understood by using the same notions as probability theory already uses, explicitly or implicitly.

Thus we come to the perhaps surprising conclusion that, in this respect, quantum systems are not as exceptional as usually thought. This thesis can be elaborated, as I started doing elsewhere (Vervoort 2011, 2012), and as I will further argue below. My position overlaps with a mathematical analysis made by L. Szabó, who comes to the conclusion that quantum probabilities are interpretable as conditional probabilities in a classical probability space (Szabó 1995, 2001).

Note that I of course do not claim that every element of quantum mechanics has its counterpart in probability theory à la von Mises. But (the observable predictions of) quantum mechanics should comply with probability theory. Probability theory has no 'uncertainty relations', neither the corresponding 'commutation relations' – two paradigmatic ingredients of quantum mechanics. Indeed, without again going in details, quantum mechanics posits that the 'commutator' of position x and momentum $p_x$ does not vanish, i.e., in symbols: $[x, p_x] = x.p_x - p_x.x = ih/2\pi \neq 0$ (here i is the imaginary number and h Planck's constant, and x and $p_x$ are operators). This is the formal basis for the well-known claim that position and momentum of a quantum particle do not exist at the same time, at least are not measurable at the same time (recall Bohr's argument against EPR). In classical mechanics, on the contrary, the commutator is zero and both properties exist at the same time.

Now, I believe that even *these commutation relations are foreshadowed by an adequate interpretation of probability à la von Mises*. As just stated, the commutation relations stipulate, among other things, which properties A and B can be measured simultaneously. On a natural reading of von Mises, P(A, B), the joint probability of A and B, is defined *only if A and B can be measured by an experimental set-up that allows to measure A, B, and A and B (combined) by using the same observing system*[31] (Vervoort 2011). This conclusion is coherent with a claim made above, namely that one can always specify the experimental conditions C for any probability and consider P(~|C) as the more precise version of P(~). According to the calculus, we have: P(A,B|C) = P(A|C).P(B|A,C). To make this a legitimate formula, C must of course represent the *same*

[31] I extract this from von Mises' exposition in his (1964) pp. 26 – 39, even if the text is not entirely explicit. The key procedure described there is the following: start from two collectives (probabilistic experimental series), one in which A is measured, one for B. If A and B are 'combinable', then one can construct a collective for the joint measurement of A and B, and determine their joint probability. On a straightforward interpretation of von Mises, the latter collective must correspond to an experimental series *using the same equipment as used for measuring A and B separately*.

experimental conditions (variables) in the three probabilities, i.e., on my model, in the three corresponding experiments defining the three probabilities in the formula. Ergo, C must refer to equipment, in general conditions of measurement/observation, allowing to measure both A and B ('simultaneously' one might say).

Thus in physical probability theory the question of simultaneous measurement of two random quantities is already crucial – just as it is crucial in quantum mechanics. Not of all quantities A and B the joint distribution exists; one can find conditions C in which A and B cannot be measured together by the same equipment. We see here again that (via the commutation relations) quantum mechanics fills in, in a stringent manner, the conditions under which to apply *classical* probability theory in the subatomic realm.

To end, I do not resist the temptation to suggest a link with the notorious 'quantum contextuality' revealed by the Kochen-Specker theorem, of which the mathematical intricacies only begin to be understood (Budroni et al. 2022). I will just be suggestive; I will just pass the word to authors who are recognised experts in this matter. Here is what Budroni et al. in their review state (2022, p. 1): "A central result in the foundations of quantum mechanics is the Kochen-Specker theorem. In short, it states that quantum mechanics is in conflict with classical models [hidden-variable theories à la Einstein] in which the result of a measurement does not depend on which other compatible measurements are jointly performed. Here, compatible measurements are those that can be performed simultaneously or in any order without disturbance. This conflict is generically called quantum contextuality." Another expert on the topic, N. D. Mermin, states[32] (1993): "The [Kochen-Specker] theorems establish that in a hidden-variables theory the values assigned even to a set of mutually commuting observables must depend on the manner in which they are measured—a fact that Bohr could have told us long ago (although he would have disapproved of the whole undertaking)."

Could von Mises have told us even earlier? Has Bohr read von Mises? These are speculations, and we will probably never know, but once more, it seems that these descriptions of quantum contextuality find their explanation in a detailed physical interpretation of probability. I believe that such an interpretation makes many quantum riddles transparent.

[32] For those interested in quantum nonlocality, the dense quote by Mermin continues thus: "And Bell's Theorem establishes that the value assigned to an observable must depend on the complete experimental arrangement under which it is measured even when two arrangements differ only far from the region in which the value is ascertained" (Mermin 1993). In this book I will not deal with the technically complex (and still utterly mysterious) notion of quantum nonlocality; but see some ideas in the Epilogue.

Note, finally, that the results **R2 – R4** are all derived by giving due attention to the experimental 'conditions' mentioned in **DEF1&2**. This is of course the main reason why they deserve to be made explicit.

## 6. Conclusion

In this chapter I have had a closer look at the decidedly subtle interpretation of probability, with the aim of proposing a unified model – a variant of von Mises' frequency interpretation. I am well aware that some topics could have been explored in more detail, and that questions will remain. Notably, it would doubtlessly be interesting to study through the lens of my model other problems that were investigated in detail by others (for instance, the case studies in Beisbart and Hartmann 2011). Unfortunately, doing so would make this chapter prohibitively long, but I believe this exercise can now well be done (and note that such case studies are not aimed at distilling a unifying framework). Next, it might be that the function 'P' of Kolmogorov's theory can be interpreted in a way that is different from mine and yet consistent. But I believe that such interpretations should have a weaker explanatory power than the model here proposed. They might be applications of the frequency interpretation valid under certain restrictive conditions – a case in point is Jaynes' Bayesianism.

To end, let us go once more over the main findings, first those that, hopefully, comply with the MCT-approach of Chapter 1, then those that are somewhat more conjectural.

### *6.1. Synthetic conclusions*

My main conclusions are that, in order to understand probability as a scientific concept, one cannot do without a frequency interpretation à la von Mises; and that probability theory, beyond the mathematics, is ultimately a physical theory, more precisely a theory about the objective world. However, it is a theory of a special kind: *it is so general that it can be abstracted into a formalised mathematical theory*. In this sense, *it seems to be the most general quantitative theory about the world we have*. It seems to me that the importance of this claim cannot be overstated; we will encounter it at several occasions in this book. In order to counter criticisms against von Mises' interpretation (which constitutes, together with Kolmogorov's formalism, the full-blown real-world theory just mentioned) I proposed a modified account based on **DEF1&2** and **C1**-**C6**. A key ingredient is that it only makes sense to define probability for *systems-in-context that exhibit frequency stabilisation* – in short: ρ-systems, experimentally repeatable in well-defined, partitioned conditions. Simply paying attention to these experimental conditions, by realising that they are a part of the ρ-system, was the key tool for solving classic problems and to come to a unified picture.

In particular, I argued that including the probing subsystem into the ρ-system allows to define probability in an objective manner, and to reframe Bayesianism à la Jaynes objectively. Similarly, the probing system of any quantum experiment, having probabilistic outcomes, is part of the ρ-system that determines the corresponding probabilities. Thus, a basic assumption of the orthodox interpretation of quantum theory finds its explanation in a more fundamental theory than quantum mechanics: (physical) probability theory. Many quantum puzzles are due to not realising that the probing system determines probability measures (on the ontic level). Realising that the initiating subsystem belongs to the ρ-system also leads to the solution of paradoxes, such as Bertrand's paradox.

Throughout this chapter we encountered a tenaciously lingering epistemic or operationalist or subjective aspect in (the application of) probability theory. The baseline is this: since a probability value can be scientifically verified by repeatedly reproducing the corresponding ρ-system (or by inferring indirect but empirically verifiable predictions), it must correspond to something objective, ontic: namely, to the objective chance of some event to occur, observed or not. Yet, we detected at several occasions an epistemic element in probability theory too, at least on the surface. Notably: *the ascription of a deterministic or indeterministic nature to a given system* does depend on the knowledge state of the attributor. But it appeared that this worry could easily be defused: the difference in ascription can be reframed as a consequence of an objective difference in experiments and ρ-systems (cf. **C2**). The conclusion seems clear: determinism or indeterminism refer not to test-systems per se, but to ρ-systems including an initialising subsystem. And a ρ-system is deterministic in an objective sense if it 'generates' only probabilities 0 and 1 (more precisely, a property X of a ρ-system is deterministic if the ρ-system generates only probabilities 0 and 1 for X). So there is also a simple criterion to decide whether a property/variable X of a system T is deterministic or not: it is only deterministic if a repeatable 'X-measuring' experiment exists on T that would always lead to trivial probabilities for X. I further argued that the counterfactual flavour of von-Mises-like definitions is a manifestation of a harmless and unavoidable aspect of any physical/real-world theory. Finally, I argued that one can avoid the notion of limit-at-infinity by introducing the notion of precision in the definition.

### *6.2. Lines of further research*

I made a few conjectures in this chapter that might be interesting to investigate further. Notably: (i) Not only Jaynes' but all coherent subjective interpretations can be unified with von Mises-style frequency models (cf. **R3** above); (ii) The study of the probabilistic foundation of quantum contextuality can be pushed forward by a mathematical analysis (cf. **R4** above); (iii) The relation between precision-intervals and causal variables determining these intervals can be pushed forward by a mathematical analysis (**C4**). Conjecture (i) is already backed-up by the explanatory power of the account presented.

Surely of a more general philosophical interest is the following matter, which concerns the universal scope of probability theory and the origin of the subjectivist temptation. Going over the main arguments one last time seems to warrant intriguing even if somewhat speculative conclusions. First, note that the claim that probability is best defined relative to a precision δ corroborates the inclusiveness of probability theory: it can also be applied to systems that cannot be repeated *very* many times in exactly the same conditions. Now, this general applicability of probability theory seems further supported by following argument. All physical systems in the narrow sense (those in physics textbooks) – or if one prefers, series of experimental outcomes obtained on such systems – are covered by probability theory (possibly restricted to trivial probability values, 0 or 1). Since we have to assume that *all* real-world systems ultimately consist of physical building blocks, it seems compelling to think that probability is applicable, *in principle*, to *all* real-world events – including events involving free-willed humans! So, is it then legitimate, after all, to speak of say the probability of 'Alice says <hi, Bob> today between 4pm and 5pm'? Recall, I have argued that probability theory only applies to ρ-random systems. And most random-looking events created by humans are certainly not ρ-random, since we do not know conditions in which these events would frequency-stabilise. But what about all-knowing experimenters? Could they identify conditions of stabilisation for any property of any system and for any event? Are there *in principle* conditions for frequency stabilisation for any event, any property, and any system, including those studied in say human and social sciences? It seems, then, that I should specify my original claim and state that probability theory *in scientific practice* only applies to ρ-random, i.e. frequency-stabilising systems, *so when verifiable numerical values are demanded*. That is why I used the notion of scientific or empirical probability, the topic of this chapter. Perhaps it makes sense, *in philosophical arguments about principles*, to use the notion of probability in a wider sense. The notion of quasi-probability that I used a few times above could then also be justified by above rationale about the universality of the concept of probability.

Once we are launched on this track, we might even go a little further. Indeed, the uniquely wide scope of probability theory is also corroborated by the fact that it can be abstracted into a purely mathematical, axiomatised and self-standing theory – as Kolmogorov has shown. This is not the usual fate of physical theories, of course; in that sense it is not exactly a physical theory in the usual meaning. It seems to me that this idea warrants a conclusion of great importance: namely, that probability theory is almost on a par with logic regarding scope of application. I believe it is one of those theories that would survive even when other basic theories, say quantum mechanics, would be replaced by others.

Now, if, as hypothesised, probability theory (ultimately, even if not practically) applies to all physical systems in a wide sense, then it should also apply to us, experimenters. But then the operational touch, the use of the phrasing "it is possible to do an experiment" in **DEF1&2** might become unavoidable, or at least transparent. I suspect one can understand it this way: probability theory also describes us, as experimenters, via the 'conditions' that play such a key role in the ontology of probability. More specifically: *these conditions prescribe how we should do 'good' experiments, which lead, by definition, to objective measures*. Only when one repeats the 'right' experimental conditions – which can be extremely sophisticated, think of quantum experiments – probabilities can be measured, and can we conclude that the measurement was well done. Deterministic experiments leading to trivial probabilities – those that occupy most of physicists' time in the laboratory – are just a special case. From this point of view, *probability theory would be the (most general) theory of scientific experimentation, of measurement.* (To the least, such a theory can immediately be derived from a von Mises-type model as ours. The key idea simply is that good measurements are operations that lead to frequency stabilisation as in **DEF1**.) This would be, at the meta-level, the reciprocal view to the subjectivist claim that probability theory is the theory of scientific inference. That probability theory is, as a mathematical theory, a 'measure-theory' is then a happy coincidence – or not.

Again, within scientific realism, this operationalist or epistemic flavour is certainly not in contradiction with the claim that probability theory tells us, at the same time, something objective about the world, independently of human experimenters; probabilities tell us something objective about systems too. But then *there should be something in nature that makes systems behave according to the rules of probability theory; there is something that stabilises frequencies*. What could this something be? (I conjecture that a coherent answer is possible; some hints are proposed in the Epilogue, Chapter 9.) Other closely related questions are: Why is frequency stabilisation a ubiquitous property of nature? Does the operational flavour of the definition of probability point to our human cognitive limitations, as Laplace thought? Is nature ultimately deterministic rather than indeterministic?

It is my hope that the notions studied in this book will help answering at least some of these fascinating questions.

### Appendix 4.1. Alleged mathematical problems of von Mises' theory

Let us have a very brief look at some other criticisms of von Mises' mathematical theory, which is based on the concept of collective. One criticism, stating that von Mises' probability as a limit for $n \to \infty$ is not always well defined (see e.g. Richter 1978) was rebutted in the main text. It becomes transparent when one realises that von Mises' theory is a physical theory, not a strictly mathematical one, as also argued in the main text (see von Mises' own arguments in his (1928/1981), e.g. p. 85). A tougher critique concerns the 'condition of randomness' that von Mises imposes on collectives. According to him, the randomness of a collective can be characterised as 'invariance under place selection': roughly, the limiting frequencies of a collective should remain invariant in subsequences obtained under 'place selections', certain functions defined on the original collective. But which and how many place selections are required ? – von Mises' critics ask. An important result was obtained by A. Wald, who showed that for any denumerable set of functions performing the subsequence selection, there exist infinitely many collectives à la von

Mises: a result that seems amply satisfying for proponents of von Mises' theory (see the reviews in von Plato (1994), Gillies (2000) and Khrennikov (2008) p. 25). A well-known objection by J. Ville can be shown not to be a real problem (see Khrennikov (2008) p. 27 and also Ville's own favourable conclusion reproduced in von Plato (1994) p. 197).

Although I believe one does not need the notion of collective and its calculus, note that von Mises' attempts to mathematically describe randomness led to interesting developments in the mathematics of string complexity, as produced by his 'competitor' Kolmogorov, and mathematicians as Martin-Löf. A general and extremely interesting result of these developments can be stated as follows: *real physical randomness defies full mathematical characterisation* (e.g. Khrennikov (2008) p. 28). This seems not really surprising; in any case, if a series of experimental results can be generated by an algorithm, it is not random by definition (it may of course look random to the naïve observer). A real series of outcomes of coin tosses cannot be generated (predicted) by an algorithm; but it does show frequency stabilisation.

This remark allows to counter a curious critique by Fine (1973), who derives a theorem (p. 93) that is interpreted by the author as showing that frequency stabilisation is nothing objective, but "the outcome of our approach to data". However, closer inspection shows that Fine's argument only applies to complex mathematical sequences (of 0's and 1's) that can be generated by computer programs. But again, if a sequence can be generated by a computer algorithm, it is by definition not random, but deterministic, even if it looks random. The essential point of real series of coin tosses is that they cannot be predicted by numerical algorithms… *and* that they show frequency stabilisation. From this perspective, mathematical randomness or complexity, as discussed in relation to the complexity of number strings, has little to do with physical randomness[33].

Note that the model I proposed is immune to criticisms against the calculus of collectives, since I do not use it. I argued that von Mises' randomness condition is not necessary for the theory, since the axioms of Kolmogorov's measure-theoretic approach can be proven without using that condition (an explicit proof in Gillies (2000) p. 112).

[33] Except if the number strings are generated by a real physical number generator, which is almost never the case.

**Appendix 4.2. Probability of the velocity of molecules and other 'theoretical' probabilities**

To further clarify **DEF1&2**, it is instructive to verify these definitions in the case of a natural probabilistic property as the velocity of gas molecules. This is a case to which von Mises devotes attention himself (von Mises 1928/1981 p. 20), since it is slightly subtle. The probability that a gas molecule has a certain velocity exists according to **DEF1&2**. Indeed, the system (the gas molecule in its environment), or more precisely the event consisting in the molecule having a certain velocity in the given environment, exhibits frequency stabilisation; according to physics it is – *at least in principle* – possible to repeat velocity measurements on the same or on similar molecules. The probabilities are for many types of gases given by the Maxwell-Boltzmann distribution, stating that the probability P($v$) that the gas molecule has a velocity in an interval $dv$ around $v$ is:

$$\mathrm{P}(v) \;=\; 4\pi . \sqrt{\left(\frac{m}{2\pi kT}\right)^{3}} \; . v^{2} . \exp\left(\frac{-mv^{2}}{2kT}\right) . dv , \qquad (2)$$

where m is the mass of the molecule, k Boltzmann's constant and T the temperature of the environment. (This implies that the probability density $\rho(v) = \mathrm{P}(v)/\mathrm{d}v$.)

This case is somewhat subtle: the probabilities (2) are known to be a good description of reality since more than a century, even if in practice they were not directly determined by velocity measurements, due to technical limitations in instrumentation. How then did physicists come to accept (2) as the correct formula ? Initially, Maxwell and Boltzmann had derived (2) based on theoretical arguments. But the essential point is that (2) could subsequently be tested in numerous manners in an *indirect* way: if one assumes (2) as the correct probability for velocities, one can use physical theory to derive various predictions for other variables (say kinetic energy) that are functionally related to velocity (kinetic energy = $mv^2/2$) – and it is these predictions which were tested so often that the starting hypothesis (2) got 'scientifically generally accepted', and indeed one of the pillars of statistical physics. This is, in essence, also von Mises' justification[34]. In sum, when we say in **DEF1** (ii) "the relative frequency of the outcome $R_j$ (j = 1,…, J) converges" we might say "the relative frequency of the outcome $R_j$ (j = 1,…, J), *or of a functionally related event*, converges". For our concern the most interesting point illustrated by this example is that the probability values (2) strongly depend on the environment; the essential parameter describing the environment is the temperature T of the gas.

[34] Here is von Mises' verdict (1928/1981 p. 20): "It is true that nobody has yet tried to measure the actual velocities of all the single molecules in a gas, and to calculate in this way the relative frequencies with which the different values occur. Instead, the physicist makes certain theoretical assumptions concerning these frequencies (or, more generally, their limiting values), and tests experimentally certain consequences, derived on the basis of these assumptions. Although the possibility of a direct determination of the probability does not exist in this case, there is nevertheless no fundamental difference between it and the other […] examples treated. The main point is that in this case too, all considerations are based on the existence of constant limiting values of relative frequencies […]".

I believe that, mutatis mutandis, a similar analysis can be made for the 'theoretical' probabilities that are considered by C. Beisbart (Beisbart and Hartmann 2011, Ch. 6). Beisbart investigates the interpretation of probabilistic properties as they are predicted by physical models, e.g. of Brownian motion and of the spatial distribution of galaxies. Without entering in detail, it seems **DEF1&2** can well be understood to also apply to these cases. For instance, a Brownian particle / system can be subject to repeated observations in constant-enough conditions; and experimental verification of stochastic properties as dwell-time or displacement radius could be done, and has been done, via the determination of frequencies as in **DEF1&2**. Somewhat more subtle are probabilistic parameters physicists use for describing spatial patterns of e.g. galaxies. If such a model describes the spatial distribution of galaxies in the probabilistic sense, via e.g. the probability $P_R(N)$ that a given radius R contains N galaxies (assuming for instance homogeneity), then it is again possible to apply the recipe of **DEF1&2**. One can pick 'many' (say n) different regions with radius R on the map, and count within each of these regions the number of galaxies. The experimental ratio "(number of regions containing N galaxies) / n" should approximate $P_R(N)$ within a workable precision if one increases n; if not physicists would not accept the model as an adequate description.

### Appendix 4.3. Illustration of the relativity of (in)deterministic ascription[35]

In many cases the relativity of (in)deterministic ascription is obvious. Consider following thought experiment, a mechanised version of a die toss, consisting of say N runs (tosses). The test involves a 20x20x20 $cm^3$ cube made of soft plastic, positioned at the beginning of each toss with high precision in the middle of a drumhead. Underneath the drumhead a metallic pin moves quickly upward, imparting to the cube a vertical momentum so that it tosses around maybe once or twice; each such toss is initiated by pressing a button. The outcome or event 'e' is the number ∈

[35] This example is extracted from a preprint (Vervoort 2020).

{1,2,…,6} on the upper face, 'measured' after landing back on the drumhead. Suppose that Alice wishes to ascertain whether this experiment is genuinely (i.e. non-trivially) probabilistic and, if so, what the values of the probabilities P(e) are; also suppose that the pin and its mechanism are screened off, so that Alice has no access to them and does not know their functioning. The best she can do is to toss the die very many times by pressing the button, and to make a table in which she notes, for each toss, the outcome e. In a trivial case she would find always the same result, say e = 6 (she is tenacious and goes until N = 10,000). In this case she would term the experiment or the system deterministic: she feels she can predict what would be the result of the (N+1)-th toss; so she assumes P(6) = 1 and P(1) = P(2) =… = 0. But consider now the case where each outcome looks perfectly random, unpredictable to her. Then she has to determine, based on her table, the frequencies $P_N(e) = \#_N(e) / N$, where $\#_N(e)$ is the number of times the outcome is e in N trials. Suppose she finds that for N = 10,000, $P_N(e)$ comes close to 1/6 for all e, and moreover that, when comparing the $P_N(e)$ for N = 100, 500, 1000 etc. the $P_N(e)$ come closer and closer to 1/6. Realising that this frequency stabilisation is, besides the unpredictability of individual outcomes, the hallmark of a (non-trivial) probabilistic system, she rather confidently concludes that the system is probabilistic and that the six probabilities P(e) are = 1/6 (to good approximation). She has done her job as a physicist.

This experiment could be deterministic in disguise: an engineer, Bob, could have made the die and the system that commands the pin in such a manner that there is a perfect functional relation between the coordinates of the pin impact and the outcome 'e'. E.g., the pin hits M = 20,000 positions on the lower surface of the die, say distributed over a square grid. (Once the first M grid positions have been probed the next runs will repeat the same sequence, so M is the periodicity of the impact points.) If the system is well-calibrated and the die movement not too chaotic (the die is of soft material and rotates not more than once or twice), it is possible to establish such a functional relation between impact position and outcome, while at the same time ensuring that the outcomes 1, 2, …, 6 arise in a random-looking sequence leading[36] to P(e) = 1/6 for all e.

So this is a system that the engineer Bob and the 'informed' or 'knowledgeable' experimenter will identify as a deterministic system – each individual toss can be predicted in advance given the pin position – while the uninformed experimenter will deem it indeterministic. This is an example of relativity or subjectivity of (in)deterministic ascription regarding physical

[36] This assumes of course that there are about M/6 grid positions that lead to each of the 6 outcomes e.

systems – relativity with respect to the knowledge state of the subject inquiring about the system. Notice again that this does not mean that probability values are subjective: both the informed Bob and the uninformed Alice will measure, when asked, the same probabilities P(e) *if they do the same experiment on the same ρ-system*. In the situation sketched, Alice and Bob have done experiments with different ρ-systems: Alice with {knob + hidden pin; die; eye}, while Bob has used the die in the conditions {knob + controlled & observed pin; die; eye}. Bob is knowledgeable because he knows or has done repeated experiments to verify that each pin position leads to a constant outcome, for instance, pin position 565 always leads to e = 1, so that $P(e = 1|pin = 565) \rightarrow 1$, etc.

Note also that there is, in principle, a way for Alice to discover the deterministic nature of the system, by noting outcomes and by making N runs with N >> M, say N = n.M (n = 2 or more), with M the periodicity of the impact points. If she does so, she will notice that also the outcomes repeat themselves with a periodicity M. She can then predict the individual outcomes with a likelihood that is proportional to n. In principle she can reach *quasi*-certainty about future outcomes, and that is all a physicist can ask for.

A less contrived example of probabilities emerging from, or reducible to, deterministic dynamics is the following. A typical numerical method in statistical mechanics is the Monte-Carlo technique, quite faithfully reproducing probabilistic features of a wide variety of real-world physical systems. Now, this numerical programming technique always uses, in practice, a pseudorandom number generator, simulating physical randomness in some variable. Such a generator is in reality deterministic: the generated numbers look randomly distributed but are actually generated by a complex function – deterministic by definition. Hence a dynamics that is fully deterministic can reproduce an enormous variety of systems of stochastic mechanics; these probabilistic systems can be understood to be deterministic under the surface.

## References

Bacciagaluppi, G. (2020). 'Unscrambling subjective and epistemic probabilities', in Hemmo and Shenker (Eds.), *Quantum, Probability, Logic: The Work and Influence of Itamar Pitowsky*, Springer, pp. 49-89

Bell, J.S. (1981). "Bertlmann's Socks and the Nature of Reality", Journal de Physique, 42, Complément C2, pp. C2-41 – C2-62

Beisbart, C. and Hartmann, S. (2011). *Probabilities in Physics* (Editors), Oxford University Press

Beisbart, C. (2011). "Probabilistic modeling in physics", in Beisbart and Hartmann (Eds.), *Probabilities in Physics*, Oxford University Press, pp. 143 – 167

Budroni C., A. Cabello, O. Gühne, M. Kleinmann, J.-Å. Larsson (2022). "Kochen-Specker contextuality", arXiv preprint arXiv:2102.13036

Bohr, N. (1935). "Quantum Mechanics and Physical Reality", Nature, 136, 1025-1026

Bub, J. (2007). "Quantum Probabilities as Degrees of Belief", Stud. in Hist. and Phil. of Mod. Phys., 38, 232-254

Bunge, M., (2008). "Bayesianism: Science or Pseudoscience?" International Review of Victimology, 15(2), 165–178

Bush, J. W. M., K. Papatryfonos, and V. Frumkin (2024). "The state of play in Hydrodynamic Quantum Analogs," In: Castro, P., Bush, J.W.M., Croca, J. (eds) *Advances in Pilot Wave Theory*, Boston Studies in the Philosophy and History of Science, Vol 344, Cham: Springer

Caves, C., C. Fuchs, and R. Schack (2002). "Quantum Probabilities as Bayesian Probabilities", Physical Review A, 65, 022305

Caves, C., C. Fuchs, and R. Schack (2007). "Subjective Probability and Quantum Certainty", Stud. in Hist. and Phil. of Mod. Phys., 38, 255-274

Couder, Y., S. Protière, E. Fort, and A. Boudaoud (2005). "Walking and orbiting droplets", Nature, 437, 208

Couder Y, and E. Fort (2006). "Single particle diffraction and interference at a macroscopic scale", Physical Review Lett., 97, 154101

Eddi, A., E. Sultan, J. Moukhtar, E. Fort, M. Rossi, and Y. Couder (2011). "Information stored in Faraday waves: the origin of a path memory", J. Fluid Mechanics, 674, 433–463

Einstein, A., B. Podolsky, and N. Rosen (1935). "Can quantum-mechanical description of physical reality be considered complete?", Physical Review, 47, 777–780

Feller, W. (1991). *An Introduction to Probability Theory and its Applications*, John Wiley & Sons

Fine, A., and T. A. Ryckman (2020). "The Einstein-Podolsky-Rosen Argument in Quantum Theory", The Stanford Encyclopedia of Philosophy (Summer 2020 Edition), Edward N. Zalta (ed.), URL = <https://plato.stanford.edu/archives/sum2020/entries/qt-epr/>

Fine, T. (1973). *Theories of Probability*, Academic Press

Gillies, D. (2000). *Philosophical Theories of Probability*, Routledge

Gnedenko, B. (1967). *Theory of Probability*, Chelsea Publishing Co.

Hájek, A. (2009). "Fifteen arguments against hypothetical frequentism", Erkenntnis 70, 211–35

Hájek, A. (2023). "Interpretations of Probability", The Stanford Encyclopedia of Philosophy (Winter 2023 Edition), Edward N. Zalta and Uri Nodelman (eds.), URL = <https://plato.stanford.edu/archives/win2023/entries/probability-interpret/>.

Jaynes, E. T. (2003). *Probability Theory. The Logic of Science*, Cambridge University Press

Jaynes, E. T. (1989). “Clearing up mysteries - the original goal”, in *Maximum Entropy and Bayesian Methods*, ed. J. Skilling, Dordrecht: Kluwer Academic, 1-27

Khrennikov, A. (2008). *Interpretations of Probability*, de Gruyter

Kolmogorov, A. (1933/1956). *Foundations of the Theory of Probability* (2nd ed.), Chelsea

Marinoff, L. (1994). “A resolution of Bertrand’s paradox”, Phil. of Science, 61, 1 – 24

Mermin, N. D. (1993). “Hidden variables and the two theorems of John Bell”, Reviews Mod. Physics 65, 803

Neapolitan, R. (1992). “A Limiting Frequency Approach to Probability Based on the Weak Law of Large Numbers”, Philosophy of Science, 59(3), 389-407

Nikolaev, V. and Vervoort, L. (2023). “Aspects of Superdeterminism Made Intuitive”, Foundations of Physics 53, 17

Papatryfonos, K., L. Vervoort, A. Nachbin, M. Labousse, and J. W. M. Bush (2024). “Static Bell test in pilot-wave hydrodynamics”, Physical Review Fluids 9(8), 084001

Peres, A. (2002). *Quantum Theory: Concepts and Methods*, Kluwer Academic Publ.

Popper, K. (1957). “The Propensity Interpretation of the Calculus of Probability, and Quantum Mechanics”, pp. 65-70, in S. Körner (ed.), Observation and Interpretation, Academic Press

Popper, K. (1959). “The Propensity Interpretation of Probability”, British Journal for the Philosophy of Science, 10, 25–42

Richter, E. (1978). *Höhere Mathematik für den Praktiker*, Johann Ambrosius Barth

Schiffrin, J.S., and R.M. Wald (2012). “Measure and probability in cosmology”, Physical Review D, 86, 023521

Smerlak, M. and Rovelli, C. (2007). “Relational EPR”, Foundations of Physics, 37, 427-445

Szabó, L. E. (1995). “Is quantum mechanics compatible with a deterministic universe ? Two interpretations of quantum probabilities”, Foundations of Physics Letters, 8, 421

Szabó, L. E. (2001). “Critical reflections on quantum probability theory”, in M. Rédei, M. Stoeltzner (Eds.), *John von Neumann and the Foundations of Quantum Physics*, Vienna Circle Institute Yearbook, Vol. 8, Kluwer

Tijms, H. (2004). *Understanding Probability: Chance Rules in Everyday Life*, Cambridge University Press

Uffink, J. (2011). “Subjective Probability and Statistical Physics” in Beisbart and Hartmann (Eds.), *Probabilities in Physics*, Oxford University Press, pp. 25 – 49

van Fraassen, B. (1980). *The Scientific Image*, Clarendon Press

Vervoort, L. (2011). “The interpretation of quantum mechanics and of probability: Identical role of the ‘observer’”, arXiv preprint https://doi.org/10.48550/arXiv.1106.3584

Vervoort, L. (2012). “The instrumentalist aspects of quantum mechanics stem from probability theory”, American Institute of Physics Conference Proceedings FPP6 (Foundations of Probability and Physics), Ed. M. D’Ariano et al., pp. 348-354

Vervoort, L. (2019). “Probability Theory as a Physical Theory Points to Superdeterminism”, Entropy, 21(9), 848 (1-13)

Vervoort, L. (2020). “The hypothesis of ‘hidden variables’ as a unifying principle in physics”, Preprint, cf. http://philsci-archive.pitt.edu/16781/

von Mises, R. (1928/1981). *Probability, Statistics and Truth*, 2nd revised English edition, Dover Publications

von Mises, R. (1964). *Mathematical Theory of Probability and Statistics*, Academic Press

von Plato, J. (1994). *Creating Modern Probability*, Cambridge University Press

Wuethrich, C. (2011). “Can the world be shown to be indeterministic after all?” In Probabilities in Physics, ed. C. Beisbart and S. Hartmann, pp. 365-389, Oxford University Press